%% file: SymplecticShadowIntegration_AIP.tex
\documentclass[%
 aip,
 amsmath,amssymb,
 reprint,%
]{revtex4-1}

\usepackage{graphicx}
\usepackage{dcolumn}
\usepackage{bm}

\usepackage{subcaption}
\usepackage{scalerel}
\usepackage{tikz}
\usetikzlibrary{svg.path}
\input{ORCID_logo.tex}
\usetikzlibrary{cd} 

\usepackage{xcolor}
\usepackage{hyperref}
\hypersetup{
	colorlinks=true,
	linkcolor={red!50!black},
	citecolor={blue!50!black},
	urlcolor={blue!20!black}
}

\usepackage{amsthm}
\usepackage{thmtools}

\usepackage[noabbrev,capitalize]{cleveref}

\usepackage[utf8]{inputenc}
\usepackage[T1]{fontenc}
\usepackage{mathptmx}
\usepackage{etoolbox}

\makeatletter
\def\@email#1#2{%
 \endgroup
 \patchcmd{\titleblock@produce}
  {\frontmatter@RRAPformat}
  {\frontmatter@RRAPformat{\produce@RRAP{*#1\href{mailto:#2}{#2}}}\frontmatter@RRAPformat}
  {}{}
}%
\makeatother

\declaretheorem[name=Remark,style=remark,numberwithin=section,qed={\hfill $\triangle$}]{remark}

\begin{document}

\def\R{\mathbb{R}}
\def\d{\mathrm{d}}
\def\p{\partial }
\def\id{\mathrm{id}}

\preprint{AIP/123-QED}

\title[Symplectic Shadow Integrators]{Symplectic integration of learned Hamiltonian systems}
\author{C. Offen $^{\protect \orcidicon{0000-0002-5940-8057}}$}%
\email{christian.offen@uni-paderborn.de}
\author{S. Ober-Blöbaum}%
\affiliation{ 
Paderborn University, Department of Mathematics, Warburger Str.\ 100, 33098 Paderborn, Germany 
}%


\date{\today}

\begin{abstract}

Hamiltonian systems are differential equations which describe systems in classical mechanics, plasma physics, and sampling problems. They exhibit many structural properties, such as a lack of attractors and the presence of conservation laws.
To predict Hamiltonian dynamics based on discrete trajectory observations, incorporation of prior knowledge about Hamiltonian structure greatly improves predictions.
This is typically done by learning the system's Hamiltonian and then integrating the Hamiltonian vector field with a symplectic integrator. For this, however, Hamiltonian data needs to be approximated based on the trajectory observations. Moreover, the numerical integrator introduces an additional discretisation error. 
In this paper, we show that an inverse modified Hamiltonian structure adapted to the geometric integrator can be learned directly from observations. A separate approximation step for the Hamiltonian data is avoided. The inverse modified data compensates for the discretisation error such that the discretisation error is eliminated. 
%
%
%
The technique is developed for Gaussian Processes.

\end{abstract}

\maketitle

\begin{quotation}
	Combining trajectory data with prior knowledge about structural properties of a dynamical system is known to greatly improve predictions of the system's motions. The article introduces the new technique {\em Symplectic Shadow Integration}, which incorporates {\em modified} structures into learned models. The modified structures compensate discretisation errors which limit the accuracy of existing approaches. 
\vspace{-0.3cm}
\end{quotation}

\vspace{-0.5cm}

\section{Introduction}

\subsection{Hamiltonian systems and symplectic integrators}

A Hamiltonian system on the phase space $M= \R^{2n}$ is a differential equation of the form
\begin{equation}
	\label{eq:HamODE}
	\dot{z} = J^{-1} \nabla H(z), \quad J=
	\begin{pmatrix}
		0 & -I_n \\ I_n & 0
	\end{pmatrix},
\end{equation}
where $H \colon M \to \R$ and $I_n$ is an $n$-dimensional identity matrix. Hamiltonian systems arise in classical mechanics, plasma physics, electrodynamics, sampling problems (Hamiltonian Monte Carlo methods) and many other applications.
Trajectories of Hamiltonian systems conserve $H$. In classical mechanics this corresponds to energy conservation. The dynamical system has no attractors. Moreover, its flow map $\phi_t \colon M \to M$ is symplectic, i.e.\ it fulfils
\begin{equation}\label{eq:SymplecticCond}
	\phi_t'(z)^\top J \phi_t'(z) = J \qquad \forall z \in M,
\end{equation}
where $\phi_t'(z)$ denotes the Jacobian matrix of $\phi_t$ at $z$. The symplectic structure has the remarkable effect that (symplectic) symmetries of $H$ yield conserved quantities of the flow by Noether's theorem.  
It is, therefore, rewarding to preserve symplecticity and symmetries when discretising \eqref{eq:HamODE} such that the numerical flow inherits qualitative features such as conserved quantities, complete integrability, and regular and chaotic regions from the exact flow. Integrators which preserve symplectic structure are called {\em symplectic integrators}. A key feature is that the numerical trajectories obtained by a symplectic integrator conserve a modified Hamiltonian or shadow Hamiltonian, which can be calculated explicitly using backward error analysis techniques.\cite{GeomIntegration}
Moreover, symplectic maps preserve phase space volume. This makes symplectic integration schemes relevant for sampling techniques such as Hamiltonian Monte Carlo methods.\cite{betancourt2018conceptual}
Next to these effects concerning initial value problems, preservation of symplectic structure is crucial to capture the bifurcation behaviour of solutions to boundary value problems in Hamiltonian systems,\cite{numericalPaper,bifurHampaper,PhDThesis} variational PDEs, \cite{PDEBifur} and in optimal control problems.\cite{offen2021bifurcation}

\subsection{Symplectic Shadow Integration (SSI)}

Techniques to identify Hamiltonian functions $H$ from data were developed by Bertalan et.\ al..\cite{Bertalan2020} For training it relies on the availability of derivatives of the flow map, so it cannot learn $H$ from observed trajectories directly.
We will modify the idea and identify an {\em inverse modified Hamiltonian} $\overline{H}$ directly from the data without any numerical approximations of derivatives.
Inverse modified equations and inverse modified Hamiltonians were introduced by Zhu et.\ al.,\cite{Zhu2020} where they are used as an analysis tool for neural networks.
An inverse modified Hamiltonian $\overline{H}$ is adapted to a symplectic integrator such that if the integrator is applied to $\dot{z} = J^{-1} \nabla \overline{H}(z)$ the discretisation error of the numerical scheme gets compensated for.

We introduce the following procedure, which we coin {\em Symplectic Shadow Integration}:

\begin{enumerate}
	\item {\em Preparation.} Choose a symplectic integrator and a step size $h$ compatible with the discrete trajectory observations.
	\item {\em Inverse system identification.} Learn the inverse modified Hamiltonian $\overline{H}$ from data.
	\item {\em Integration.} Apply the symplectic integrator to the inverse modified Hamiltonian system $\dot{z} = J^{-1} \nabla \overline{H}(z)$ to obtain a numerical flow map.
	\item {\em System identification.} If required, compute $H$ from $\overline{H}$ using backward error analysis techniques to obtain physical insight or for verification.
	
\end{enumerate}

\subsection{Relation to other approaches}\label{sec:advantages}

Before we review the idea of backward error analysis and provide details on the steps of the Symplectic Shadow Integration (SSI) procedure, let us outline the advantages of the integration technique over other more direct approaches to predict Hamiltonian dynamics from data and contrast SSI to techniques in the literature. 

\subsubsection*{Comparison to learning the flow map directly (Strategy 1).}

One could use established learning techniques, such as artificial neural networks, Gaussian processes, or kernel methods to learn the flow map of the system directly from trajectory data. SSI has the following advantages over this approach.

\begin{itemize}
	\item
	Hamiltonian structure is incorporated into the learned system. This guarantees important qualitative aspects of the prediction such as energy conservation, preservation of phase space volume and topological properties of the phase portrait.
	
	\item
	Only a real valued map $\overline{H}$ needs to be learned rather than the flow map, reducing the dimension of the learning problem and data requirements.
	
	\item
	Hamiltonian structure can be identified. It provides physical insight into the dynamics and can be used for verification. More precisely, the predicted motions using the SSI technique are the exact motions of an identified Hamiltonian, which can be computed from $\overline{H}$. The performance of SSI can, therefore, be evaluated using backward error analysis.
	
	\item
	SSI provides a framework to incorporate further prior knowledge about conservation laws, such as (angular) momentum conservation, through a combination with symmetric learning, for instance using symmetric kernels for Gaussian processes.\cite{ridderbusch2021learning} The conservation laws are then guaranteed by a discrete Noether theorem.
	
\end{itemize}

In this context, we mention Symplectic Neural Networks (SympNets).\cite{SympNets} SympNets can be used to learn the flow map of Hamiltonian systems and the learned map is guaranteed to be symplectic. In contrast, SSI learns a scalar valued map $\overline{H}$ related to the Hamiltonian of the system. Next to neural networks, SSI can be used with Gaussian Processes and kernel methods, which will be the focus of this work. A technique analogous to SSI has been developed for artificial neural networks in the recent preprint.\cite{david2021symplectic}

Another technique to incorporate geometric structure is to learn a generating function of the symplectic flow map.\cite{Rath2021,Toa2021} The learned flow map is then guaranteed to be symplectic. In contrast, SSI identifies Hamiltonian structure, albeit the approaches can yield the same learning problem for some simple numerical schemes. However, SSI uses higher order terms to correct predictions of the Hamiltonian function in a post processing step.

\subsubsection*{Comparison to learning the exact Hamiltonian and then using a symplectic integrator (Strategy 2).}

Techniques have been developed to learn the Hamiltonian rather than the flow map from data using Gaussian process regression or artificial neural networks.\cite{Bertalan2020}
This approach requires a subsequent application of a classical numerical integrator to predict motions.
SSI has the following advantages over learning the exact Hamiltonian $H$ and applying a symplectic integrator to \eqref{eq:HamODE}.

\begin{itemize}
	
	\item
	The numerical integrator introduces a discretisation error in addition to uncertainty in the Hamiltonian due to limited training data. SSI compensates this discretisation error such that high accuracy and excellent energy behaviour can be achieved despite large step sizes.
	
	\item
	Step size selection is decoupled from accuracy requirements. This is beneficial if the learned Hamiltonian and its gradient are expensive to evaluate.
	
	\item
	SSI uses the trajectory data directly. There is no need to approximate data of the underlying vector field, which would include additional discretisation errors.
	
\end{itemize}

Hamiltonian neural networks \cite{NEURIPS2019_26cd8eca} constitute an example for a technique covered by strategy 2: the Hamiltonian of the dynamical system is parametrised as an artificial neural network. For training, observations of velocities and derivatives of conjugate momenta are required. If a special form of the Hamiltonian is assumed, such as the mechanical form $H(q,p) = p^\top M(q) p + V(q)$, the components $M$ and $V$ can be parametrised separately by neural networks.\cite{zhong2020symplectic} Exploiting the mechanical form of the Hamiltonian, only position and velocity observations are required. Moreover, control terms can be added. In contrast, SSI applies to all canonical Hamiltonian systems, does not assume a particular form of the Hamiltonian, and requires observations of position and momentum data only. Observations of derivatives of these quantities are not needed.

Hamiltonian neural networks can suffer from two types of discretisation errors: the training data needs to contain velocity data and information on the derivative of conjugate momenta. Typically, these need to be approximated from observed position and momentum data, which introduces a discretisation error. Another discretisation error occurs when the learned Hamiltonian is integrated using a numerical method. These errors have been analysed by Zhu et.\ al..\cite{Zhu2020} SSI can be trained directly on position and momentum data and all discretisation errors are compensated, which is important when high accuracy is required.

A different approach to compensate discretisation errors can be found in the work by Park et.\ al..\cite{Park_2020} It applies to spin systems: a neural network is trained and used to compensate the discretisation error of a numerical time stepping scheme. Therefore, large time steps can be used and computations can be accelerated. In contrast, SSI uses analytical methods to correct discretisation errors and applies to all canonical Hamiltonian systems. Again, large time steps can be used without loosing accuracy.
Moreover, SSI does not require explicit knowledge of the Hamiltonian but identifies Hamiltonian structure from data.



The paper proceeds as follows: \cref{sec:BEA} contains a review of (inverse) modified Hamiltonians and of backward error analysis techniques. After a brief review of Gaussian process regression, \cref{sec:learnH} explains how to learn inverse modified Hamiltonians using Gaussian processes for SSI. The paper proceeds with numerical experiments in \cref{sec:Numerical} and concludes with remarks on possible extensions of SSI in \cref{sec:FutureWork}.

\section{(Inverse-) Backward error analysis (BEA)}\label{sec:BEA}

Backward error analysis (BEA) is a well-established tool of traditional numerical analysis.
Its role in this work is threefold: it provides a theoretical justification of the SSI technique, is integral when SSI is used to perform system identification, and is employed to evaluate the quality of numerical experiments.

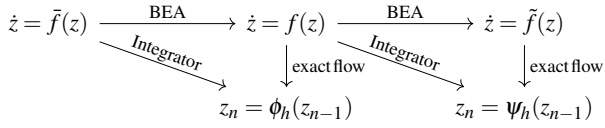
\begin{figure}
\begin{tikzcd}	[column sep={10em,between origins}]	
		\dot{z}=\bar{f}(z) \arrow{r}{\mathrm{BEA}} \arrow{rd}[anchor=center,rotate=-20,yshift=1ex]{\mathrm{Integrator}} & \dot{z}=f(z) \arrow{r}{\mathrm{BEA}} \arrow{d}{\mathrm{exact\, flow}} \arrow{rd}[anchor=center,rotate=-20,yshift=1ex]{\mathrm{Integrator}} & \dot{z}= \tilde{f}(z) \arrow{d}{\mathrm{exact\, flow}} \\
		& z_n = \phi_h(z_{n-1})                               & z_n = \psi_h(z_{n-1})
	\end{tikzcd}
	\caption{Illustration of the relation of the vector field $f$, the modified vector field $\tilde{f}$, and the inverse modified vector field $\bar{f}$.
		Backward error analysis (BEA) provides a tool to calculate the effect of an integrator on a vector field, i.e.\ to compute $\tilde{f}$ from $f$ as formal power series in the integrator's step size.
	}\label{fig:InvModf}
\end{figure}

When a differential equation $\dot{z}=f(z)$ is discretised using a numerical method, one obtains a flow map $\tilde{\phi}$ that approximates the exact flow $\phi$ of the system. The forward error measures the difference between $\phi$ and $\tilde{\phi}$. In contrast, backward error analysis seeks a modified differential equation $\dot{z}=\tilde{f}(z)$ whose exact flow map coincides with $\tilde{\phi}$. One can then compare the exact vector field $f$ with the modified (or shadow) vector field $\tilde{f}$ and understand the properties of the numerical method through the properties of $\tilde{f}$. 
The idea has been successfully applied to analyse long-term behaviour of numerical methods, especially symplectic integrators.\cite{GeomIntegration,Leimkuhler2005}
For a given integration method, an {\em inverse modified differential equation}\cite{Zhu2020} $\dot{z} = \bar{f}(z)$ is a differential equation such that the integration scheme applied to $\bar{f}$ yields the exact flow $\phi$, in other words $\tilde{\bar{f}} = f$.
The relations of $f$, $\tilde{f}$, and $\bar{f}$ are illustrated in \cref{fig:InvModf}. 

\begin{figure}
	\begin{tikzcd}[column sep={10em,between origins}]	
		\dot{z}=J^{-1} \nabla \overline{H}(z) \arrow{r}{\mathrm{BEA}} \arrow{rd}[anchor=center,rotate=-20,yshift=1ex]{\mathrm{Integrator}} & \dot{z}=J^{-1} \nabla {H}(z) \arrow{r}{\mathrm{BEA}} \arrow{d}{\mathrm{exact\, flow}} \arrow{rd}[anchor=center,rotate=-20,yshift=1ex]{\mathrm{Integrator}} & \dot{z}=J^{-1} \nabla \widetilde{H}(z) \arrow{d}{\mathrm{exact\, flow}} \\
		& z_n = \phi_h(z_{n-1})                               & z_n = \psi_h(z_{n-1})
	\end{tikzcd}
	\caption{Illustration of modified and inverse modified Hamiltonian equations. Formal power series $\overline{H}$ and $\widetilde{H}$ in the step size of the integrator exist such that the diagram commutes if the integration scheme is symplectic. }\label{fig:InvModIllu}
\end{figure}
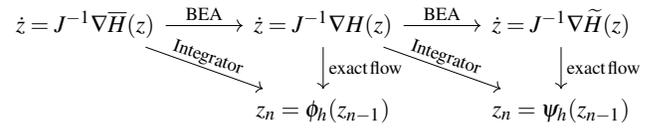

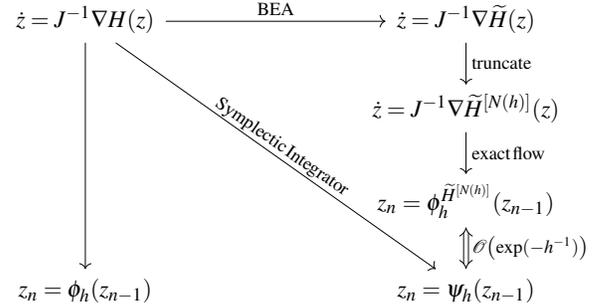
\begin{figure}
	\[
	\begin{tikzcd}[column sep={16em,between origins}]
		\dot{z}=J^{-1} \nabla H(z) \arrow{rddd}[anchor=center,rotate=-35,yshift=1ex]{\mathrm{Symplectic\, Integrator}} \arrow{r}{\mathrm{BEA}}  \arrow[ddd] & \dot{z}=J^{-1} \nabla \widetilde{H}(z) \arrow{d}{\mathrm{truncate}}   \\
		& \dot{z}=J^{-1} \nabla \widetilde{H}^{[N(h)]}(z) \arrow{d}{\mathrm{exact\, flow}}   \\
		& z_n = \phi^{\widetilde{H}^{[N(h)]}}_h(z_{n-1}) \arrow[Leftrightarrow]{d}{\mathcal{O}\left(\exp(-h^{-1})\right)} \\
		z_n = \phi_h(z_{n-1})                       & z_n = \psi_h(z_{n-1})    
	\end{tikzcd}
	\]
	\caption{\Cref{fig:InvModIllu,fig:InvModf} only show a formal analysis. The power series $\widetilde{H}$ typically does not converge. For a rigorous analysis, $\widetilde{H}$ is optimally truncated at an index $N(h)$.
		The exact flow $\phi^{\widetilde{H}^{[N(h)]}}$ of the truncated Hamiltonian $\widetilde{H}^{[N(h)]}$ agrees with the numerical flow $\psi_h$ up to an exponentially small error. }\label{fig:RigorousBEA}
\end{figure}

Symplectic integrators have the remarkable property that applied to a Hamiltonian vector field $f=J^{-1} \nabla H$ the modified vector field $\tilde{f}$ is of the form\cite{GeomIntegration} $\tilde{f}=J^{-1} \nabla \widetilde{H}$ and the inverse modified vector field $\bar{f}$ is of the form\cite{Zhu2020} $\bar{f}=J^{-1} \nabla \overline{H}$. See \cref{fig:InvModIllu} for an illustration.
The Hamiltonian structure of $\tilde{f}$ and $\bar{f}$ provides an explanation for the excellent preservation properties of symplectic integrators: the energy error $H(z_j)-H(z_0)$ oscillates within a band of width $\mathcal{O}(h^{p})$ on exponentially long time intervals if the integrator is of order $p$ and the trajectory stays within a compact set of the phase space. This needs to be contrasted with the generic energy error behaviour of non-symplectic methods which is of order $\mathcal{O}(t h^{p})$. 
%
Modified and inverse modified vector fields can be computed as formal power series in the step size $h$. Although the formal power series $\tilde{f}$ and $\widetilde{H}$ typically do not converge, they still govern the dynamics as optimal truncation results are available.\cite{GeomIntegration} The situation is illustrated in \cref{fig:RigorousBEA}.

Let us review the Symplectic Euler method and the implicit midpoint method and their traditional backward error analysis formulas. These will not only be employed to analyse the quality of our numerical schemes but will also play a role when we identify the energy of dynamical systems using SSI. Moreover, we will relate the inverse modified Hamiltonians of the two schemes with generating functions of appropriate types.

The Symplectic Euler method is a 1st order symplectic integrator. Applied to a partitioned system of differential equations
\begin{equation}
	\begin{pmatrix}
		\dot{q} \\ \dot{p}
	\end{pmatrix}
	= f (q,p), \qquad q \in \R^n, p \in \R^n
\end{equation}
with step size $h$ it reads 
\begin{equation}\label{eq:symplEuler}
	\begin{pmatrix}
		\bar{q} \\ \bar{p}
	\end{pmatrix}
	= \begin{pmatrix}
		q \\ p
	\end{pmatrix}
	+h f (\bar{q},p).
\end{equation}
Of the $2n$ equations, the first $n$ equations define $\bar{q}$ implicitly, while $\bar{p}$ can be computed explicitly from the last $n$ equations once $\bar{q}$ is known.
If $f(q,p)=J^{-1}\nabla H(q,p)$, then the first terms of $\widetilde{H}$ are given by\cite{GeomIntegration}
\begin{align}
	\label{eq:tildeH}
	\widetilde{H}
&= H
-\frac{h}{2}H_q^\top H_p\\ \nonumber
&+ \frac{h^2}{12}
\left( H_q^\top H_{pp} H_q 
+ H_p^\top H_{qq} H_p
+ 4 (H_p^\top H_{qp} H_q)\right)
+ \ldots.
\end{align}
Here $H_q$ denotes the gradient $H_q=\left(\frac{\p H}{\p q^j}\right)_{j=1}^n$ and $H_{qq}=\left(\frac{\p ^2 H}{\p q^i \p q^j}\right)_{i,j=1}^n$. $H_p$, $H_{pp}$ and $H_{qp}$ are defined analogously.
The inverse modified Hamiltonian is the Hamiltonian $\overline{H}$ such that $\widetilde{\overline{H}} = H$. Its power series is easily found by plugging the ansatz $\overline{H} = H+h \overline{H}^1 + h^2 \overline{H}^2 + \mathcal{O}(h^3)$ into $\widetilde{\overline{H}} = H$ and comparing powers of $h$. For the Symplectic Euler method we obtain
\[
\begin{split}
\overline{H} &= H + \frac{h}{2} H_q^\top H_p\\
&+ \frac{h^2}{6}
\left( H_q^\top H_{pp} H_q 
+ H_p^\top H_{qq} H_p
+ H_p^\top H_{qp} H_q\right)
+ \mathcal{O}(h^3).
\end{split}
\]
The power series $\overline{H}$, $H$, $\widetilde{H}$ are related as shown in \cref{fig:InvModIllu}.

The Symplectic Euler method with step size $h$ applied to $\overline{f}(q,p)=J^{-1} \nabla \overline{H}(q,p)$ reads
\begin{align*}
\overline{q} &= q + h \overline{H}_{\overline{p}}(q,\overline{p})\\
\overline{p} &= p - h \overline{H}_{q}(q,\overline{p}).
\end{align*}
We observe that $S_2(q,P)=q^\top P+h \overline{H}(q,P)$ is a generating function of type 2 for the exact flow map $(\overline{q},\overline{p})=\phi_h(q,p)$. Since $\phi_h$ is close to the identity, the existence of $\overline{H}$ is guaranteed\cite{GeomIntegration} for sufficiently small $h$.

Another example of a symplectic integration method is the implicit midpoint rule, which is second order accurate. Applied to the differential equation $\dot{z} = f(z)$ it reads
\begin{equation}
	\label{eq:MD}
	\bar{z} = z + h f\left(\frac{\bar{z}+z}{2}\right).
\end{equation}
If $f(z)=J^{-1}\nabla H(z)$ then the modified Hamiltonian is given as\cite{GeomIntegration}
\begin{equation}
	\label{eq:MDmodH}
	\widetilde{H} = H-\frac{h^2}{24} f^\top \mathrm{Hess}(H) f +  \mathcal{O}(h^4)
\end{equation}
and the inverse modified Hamiltonian as
\begin{equation}
	\label{eq:MDinvmodH}
	\overline{H} = H+\frac{h^2}{24} f^\top \mathrm{Hess}(H) f +  \mathcal{O}(h^4).
\end{equation}
Here $\mathrm{Hess}(H)$ denotes the Hessian matrix of $H$. Only even powers of $h$ occur in $\widetilde{H}$ and $\overline{H}$.

As for the Symplectic Euler method, the inverse modified Hamiltonian to the implicit midpoint rule can be interpreted as a generating function expressed in the coordinates $\frac{\overline{q}+q}{2},\frac{\overline{p}+p}{2}$ for the exact flow map $\phi_h$. Its existence is guaranteed for sufficiently small $h$. Details can be found in \cref{App:GenFuncMP}.


\begin{remark}
	B-series and P-series methods constitute a large class of numerical schemes which include Runge--Kutta methods and partitioned Runge--Kutta methods. Explicit formulas for inverse modified Hamiltonians for symplectic B- and P-series methods have been calculated and their expressions can be given using the theory of rooted trees (in the sense of graph theory)\cite{GeomIntegration}. Using the relation $\widetilde{\overline{H}} = H$, we derive a recursion for the terms $\overline{H} = \overline{H}_1+h\overline{H}_2+h^2\overline{H}_3+\ldots$ for consistent symplectic $B$- and $P$-series methods: in the notation of Theorem IX.9.8 in the book by Hairer et.\ al.\cite{GeomIntegration} for a $B$ series method we obtain 
	\[
	\overline{H}_1 = H, \qquad
	\overline{H}_{k+1}= \sum_{j=2}^{k+1} \sum_{\tau \in T^\ast, |\tau|=j}^{} \frac{b(\tau)}{\sigma(\tau)} \overline{H}_{k-j+2}(\tau)
	\]
	for $k=0,1,\ldots$.
	A corresponding formula for $P$ series methods can be obtained
	analogously by inverting the traditional backward error analysis formula given in Theorem IX.10.9.\cite{GeomIntegration}
\end{remark}

Applying the argumentation of traditional backward error analysis\cite{GeomIntegration} to inverse modified equations, the existence of $\overline{H}$ as a formal power series for any symplectic method was proved in the work by Zhu et. al..\cite{Zhu2020}
For non-symplectic methods, $\overline{H}$ does not exist.

In this paper, first the inverse modified Hamiltonian $\overline{H}$ will be learned from data and then the vector field $J^{-1} \nabla \overline{H}$ will be integrated using the corresponding symplectic integrator. System identification is then performed by computing $\widetilde{\overline{H}}$ to obtain $H$, i.e.\ by applying the traditional backward error analysis formulas \eqref{eq:tildeH} or \eqref{eq:MDmodH} to $\overline{H}$.

\begin{remark}
	Backward error analysis is known to describe the behaviour of numerical solutions well, not only for small but also for moderate to large time steps.\cite{GeomIntegration} We can expect Symplectic Shadow Integration (SSI) to work well for those time step sizes for which backward error analysis techniques apply. In particular, we will be able to use discretisation parameters $h$ of moderate size in the following numerical experiments.
\end{remark}



\section{Learning inverse modified Hamiltonians from trajectory data}\label{sec:learnH}

With its existence established (as formal power series in the general case and explicitly proved for the Symplectic Euler method and the implicit midpoint rule by relating their inverse modified Hamiltonians to certain generating functions), we proceed to learning an inverse modified Hamiltonian $\overline{H}$ from data using Gaussian process regression.

For an introduction to Gaussian process regression for machine learning see the book by Rasmussen and Williams\cite{Rasmussen2005}.
Let us briefly recall some relevant aspects: in Gaussian process regression a function $S \colon M \to \R$ is modelled as a sequence of random variables $(\hat{S}(z))_{z \in M}$ over an index set $M$. Each random vector $(\hat{S}(z))_{z \in \mathring{M}}$, with $\mathring{M}$ a finite subset of $M$, is multivariate normally distributed with covariance matrix $K(Z,Z) := (k(z,w))_{z,w \in \mathring{M}}$ for a covariance or kernel function $k \colon M \times M \to \R$ and mean $m(Z):=(m(z))_{z \in \mathring{M}}$ for a function $m \colon M \to \R$. For given finite data $(Z',S(Z')):=\{(z',S(z'))\}_{z' \in M'}$, mean function $m$, and kernel $k$, the posterior distribution of the random vector $\hat{S}(Z)|(Z',S(Z')) := (\hat{S}(z))_{z \in \mathring{M}} | (Z',S(Z'))$ is again multivariate normally distributed. Its mean is given as
\[
\begin{split}
\mathbb{E}&[\hat{S}(Z)|(Z',S(Z'))]\\ &= m(Z) + K(Z,Z')K(Z',Z')^{-1}(S(Z')-m(Z')),
\end{split}
\]
where $K(Z,Z')=(k(z,z'))_{z \in \mathring{M}, z' \in M'}$. The covariance matrix is given as
\[
\begin{split}
\mathrm{Cov}&[\hat{S}(Z)|(Z',S(Z'))] \\ &= K(Z,Z)-K(Z,Z')K(Z',Z')^{-1}K(Z',Z).
\end{split}
\]
The mean $\mathbb{E}[\hat{S}(Z)|(Z',S(Z'))]$ can be used as a prediction of $(S(z))_{z \in \mathring{M}}$, while the variance at a given point $z$ can be interpreted as a measure of the model uncertainty in the prediction of $S(z)$.

A technique to learn the Hamiltonian $H$ of a system from values of the Hamiltonian vector field using Gaussian processes was introduced by Bertalan et.\ al..\cite{Bertalan2020} We modify the idea such that we can use data points $(y,\phi_h(y))$ of the Hamiltonian flow $\phi$ rather than of the Hamiltonian vector field. Therefore, our method applies when there are trajectory observations available but the underlying vector field is unknown. Moreover, instead of learning $H$ we learn the inverse modified Hamiltonian $\overline{H}$. The exact Hamiltonian $H$ is then computed in a post-processing step. The following technique directly extends to all kernel methods with sufficiently smooth kernels. 

Let $Z=(z_1,\ldots,z_N)$ be $N$ points in the phase space $M$ for which we would like to predict the values $\overline{H}(Z) = (\overline{H}(z_1),\ldots,\overline{H}(z_N))$ of the inverse modified Hamiltonian corresponding to the Symplectic Euler method with step size $h$.
The prediction is based on observed data of the flow which maps a collection of $\tilde N$ points of the phase space $Y=(y_1,\ldots,y_{\tilde N}) = ((q_1,p_1),\ldots,(q_{\tilde N},p_{\tilde N}))$ to the collection of points $\bar Y=(\bar y_1,\ldots,\bar  y_{\tilde N}) = ((\bar q_1,\bar p_1),\ldots,(\bar q_{\tilde N},\bar p_{\tilde N}))$ after time $h$.
To apply Gaussian Process regression, the $2n$-dimensional phase space $M$ is interpreted as an index set. 
To derive the method, first we assume that we already have the corresponding values $\overline{H}(Z)$. This will help us to derive a linear system of equations for $\overline{H}(Z)$.

Let $\hat {\overline{H}}$ be a Gaussian process with index set $M$, a continuously differentiable kernel $k\colon M \times M \to \R$, and a constant zero mean function $m \equiv 0$. Now $\overline{H}$ can be predicted at a new point $y\in M$ as the conditional expectation 
\begin{equation} \label{eq:expectH}
	\mathbb{E}[\hat {\overline{H}}(y) | (Z,\overline{H}(Z)) ] = k(y,Z)^\top k(Z,Z)^{-1}\overline{H}(Z).
\end{equation}
Here $k(Z,Z) \in \R^{N \times N}$ is the covariance matrix $k(Z,Z)_{i,j} = k(z_i,z_j)$ and $k(y,Z) \in \R^{N \times 1}$ is given as $k(y,Z) = k(y,z_j)$. Differentiation of \eqref{eq:expectH} with respect to $y$ yields
\begin{equation} \label{eq:expectDH}
	\mathbb{E}[ \nabla \hat {\overline{H}}(y) | (Z,\overline{H}(Z)) ] = \nabla_1 k(y,Z)^\top k(Z,Z)^{-1}\overline{H}(Z) =: J \bar{f}(y)
\end{equation}
with $\nabla_1 k(y,Z) \in \R^{N \times 2n}$ given as $(\nabla_1 k(y,Z))_{i,j} = \frac{\p k}{\p y^j}k(y,z_i)$. Here $J$ is the symplectic structure matrix from \eqref{eq:HamODE}.

Let $Y=(y_1,\ldots,y_{\tilde N}) = ((q_1,p_1),\ldots,(q_{\tilde N},p_{\tilde N}))$ be a collection of $\tilde N$ points of the phase space and let $\bar Y=(\bar y_1,\ldots,\bar  y_{\tilde N}) = ((\bar q_1,\bar p_1),\ldots,(\bar q_{\tilde N},\bar p_{\tilde N}))$ denote the corresponding values of the Hamiltonian flow after time $h$ (training data). Imposing that $\bar{Y}$ was obtained from $Y$ using the Symplectic Euler method on $\bar{f}$ yields the relations
\begin{align}
\label{eq:symplELfbar}
		\begin{pmatrix}
			\bar q_j \\ \bar p_j
		\end{pmatrix}
		&= 	\begin{pmatrix}
			q_j \\ p_j
		\end{pmatrix}
		+h \bar{f}(\bar q_j, p_j)\\ \nonumber
		&=
		\begin{pmatrix}
			q_j \\ p_j
		\end{pmatrix}
		+h J^{-1} \nabla_1 k( (\bar q_j, p_j) ,Z)^\top k(Z,Z)^{-1}\overline{H}(Z)
\end{align}
for $j=1,\ldots,\tilde{N}$. 
The Hamiltonian of a Hamiltonian system is defined up to an additive constant. We can, therefore, impose $\mathbb{E}\left[ \hat {\overline{H}}(y_0) | (Z,\overline{H}(Z)) \right] = \overline{H}_0$ for any $\overline{H}_0\in \R$ (normalisation). Together with \eqref{eq:symplELfbar} we obtain a linear system
\begin{equation}\label{eq:SystemH}
	\begin{pmatrix}
		\nabla_1 k( (\bar q_1, p_1) ,Z)^\top k(Z,Z)^{-1}\\
		\vdots \qquad \qquad \vdots \\
		\nabla_1 k( (\bar q_{\tilde N}, p_{\tilde N}) ,Z)^\top k(Z,Z)^{-1} \\
		k(y_0,Z)^\top k(Z,Z)^{-1}
	\end{pmatrix}
	\begin{pmatrix}
		\overline{H}(z_1)\\ \vdots \\ \overline{H}(z_N)
	\end{pmatrix}
	= \frac{1}{h}\begin{pmatrix}
		J(\bar{y}_1 - y_1) \\ \vdots \\ 	J(\bar{y}_{\tilde N} - y_{\tilde N}) \\ \overline{H}_0
	\end{pmatrix}
\end{equation}
for $\overline{H}(Z) = \begin{pmatrix}
	\overline{H}(z_1),& \ldots, & \overline{H}(z_N)
\end{pmatrix}^\top$, which consists of $n \tilde{N}+1$ equations for $N$ unknowns. The last equation corresponds to the normalisation.

We conclude that given data $Y$,$\bar{Y}$ of the flow map of a dynamical system, e.g.\ obtained from observations of trajectories at times $\tau$ and $\tau+h$, values of the inverse modified Hamiltonian $\overline{H}(Z)$ can be predicted over points $Z$ of the phase space by solving \ref{eq:SystemH} in the least square sense.
The solution has the following interpretation: 
if \eqref{eq:SystemH} is solved exactly, the collection $(Z,\overline{H}(Z))$ has the property that if the mean of the Gaussian process $\hat{\overline{H}} | Z,\overline{H}(Z)$ is used to predict $\overline{H}(Y)$ then an application of the Symplectic Euler method to $J^{-1}\nabla\overline{H}(y_j)$ at $y_j$ recovers $\bar{y}_j$ for all $j=1,\ldots,\tilde{N}$.

\section{Numerical experiments}\label{sec:Numerical}

We apply Symplectic Shadow Integration (SSI) with the Symplectic Euler method (SE) and midpoint rule (MP) to the mathematical pendulum and the H\'enon--Heiles system. To learn the inverse modified Hamiltonian, we employ Gaussian Process regression with radial basis functions
\begin{equation}\label{eq:RBF}
	k(x,y) = k_c \exp \left( -\frac{1}{e^2} \|x-y\|^2 \right)
\end{equation}
as kernels.
When analysing the quality of our numerical results, we particularly focus on whether the phase portrait topology has been captured correctly such that long-term predictions yield qualitatively correct results. This is done by applying backward error analysis to compute the system whose exact flow coincides (up to a truncation error) with the SSI prediction, which is then compared to the exact system. 

For comparison to other approaches, we fit a Gaussian process (GP) directly to the training data $(Y,\bar{Y})$ using again radial basis functions as kernels, where the parameters $k_c, e$ are fitted using marginal likelihood estimation. For this, we employ the Python package {\tt scikit-learn}.\cite{ScikitLearn} This corresponds to strategy 1 of \cref{sec:advantages} (learning the flow map directly without using Hamiltonian structure).
For comparison with strategy 2 type approaches, we apply the Symplectic Euler method or midpoint rule to the exact system with the same step size as used in the SSI scheme. This corresponds to strategy 2 of \cref{sec:advantages} in the infinite data limit, where the Hamiltonian has been learned up to machine precision.

Source code can be found in our GitHub repository.\cite{ShadowIntegratorSoftware}

\subsection{Mathematical Pendulum}\label{sec:Pendulum}

We consider the Hamiltonian system \eqref{eq:HamODE} of the mathematical pendulum with Hamiltonian
\[
H(q,p)=\frac{1}{2}q^2 + (1-\cos(q)).
\]
To generate training data, we use a Halton sequence $Y=(y_1,\ldots,y_N)$ on the phase space $M = [-2\pi,2\pi] \times [-1.2,1.2]$ of length $N=160$. Values $\bar{Y}=(\bar{y}_1,\ldots, \bar{y}_N)$ are obtained by integrating \eqref{eq:HamODE} up to time $h=0.3$ with high precision. For this we use $n_{\mathrm{LP}}=800$ steps of the St{\"o}rmer-Verlet scheme.\cite{GeomIntegration}
As a kernel for the GP we use radial basis functions \eqref{eq:RBF} with parameters $k_c=1$ and $e=2$.
Setting $Z=Y$ in \eqref{eq:SystemH} and the step size $h=0.3$ we obtain values $\overline{H}(Z)$. We can now compute values for $\overline{H}(z)$ or $\nabla \overline{H}(z)$ for any $z \in M$ using \eqref{eq:expectH} or \eqref{eq:expectDH}, respectively.
For numerical stability, a Tikhonov regularisation $k(Z,Z) \mapsto k(Z,Z) + \sigma I_{nN}$ with $\sigma = 10^{-13}$ is applied and {\tt SciPy}'s\cite{SciPy} Cholesky solver is used, wherever multiplication with the inverse covariance matrix $k(Z,Z)^{-1}$ is required.
Higher derivatives of $\overline{H}$ are obtained using derivatives of \eqref{eq:expectDH}.

\begin{figure}

	\subcaptionbox{Phase plot of trajectories }{
		\includegraphics[width=0.4\linewidth]{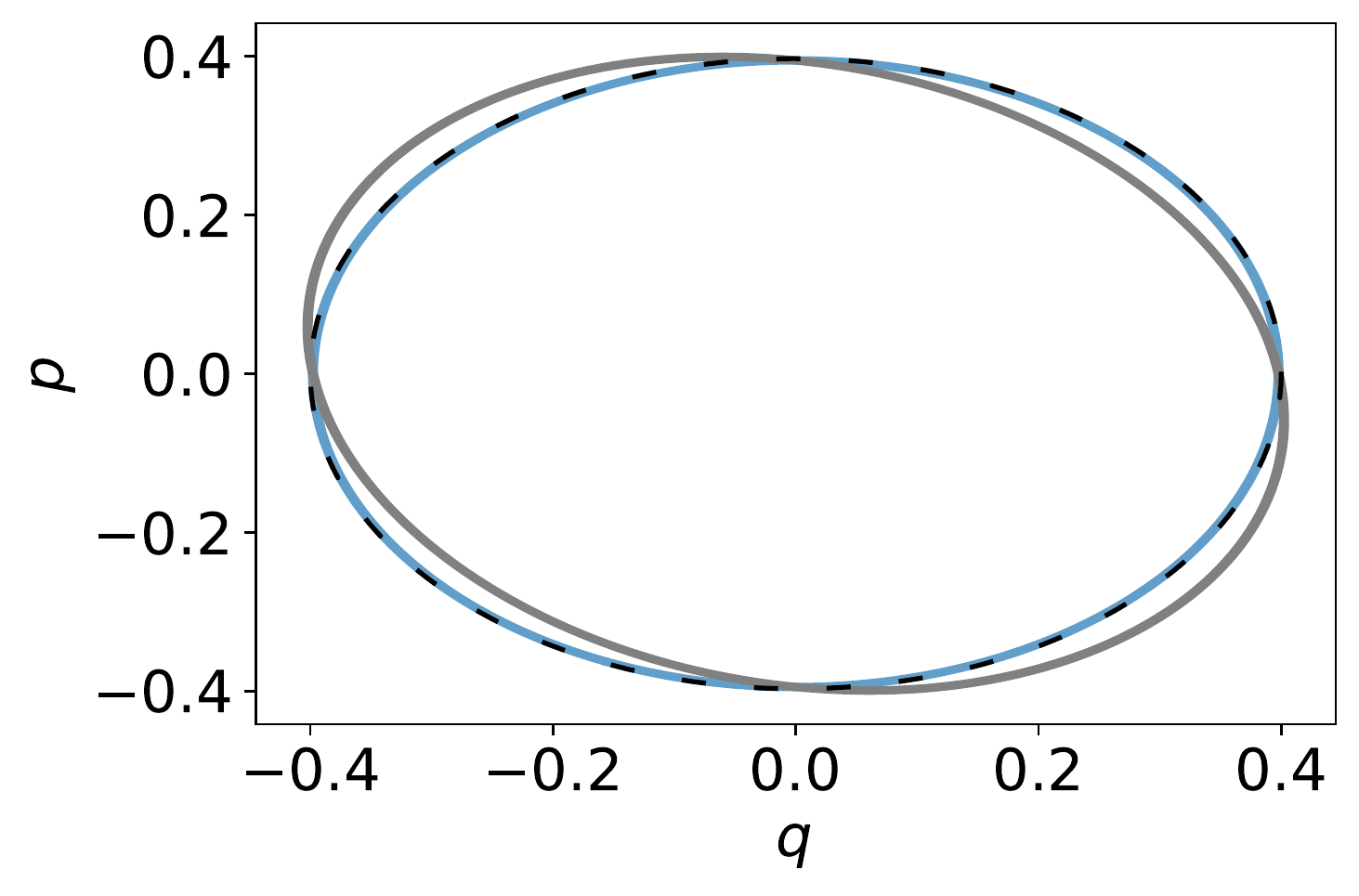}
	}
	\subcaptionbox{Conservation of $H$}{
		\includegraphics[width=0.44\linewidth]{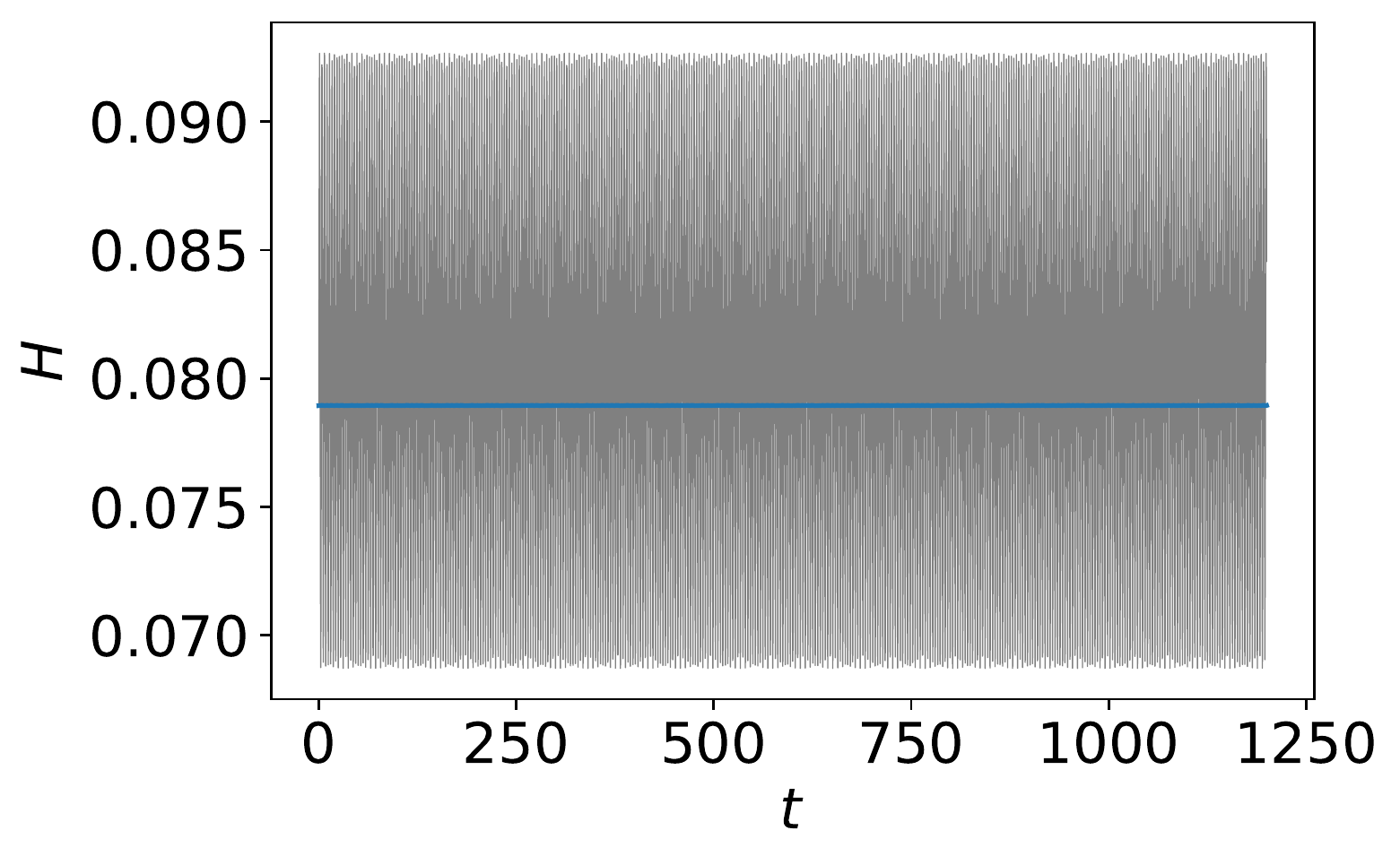}
	}\\
	\subcaptionbox{Conservation of $H-\mathrm{mean}(H)$}{
		\includegraphics[width=0.4\linewidth]{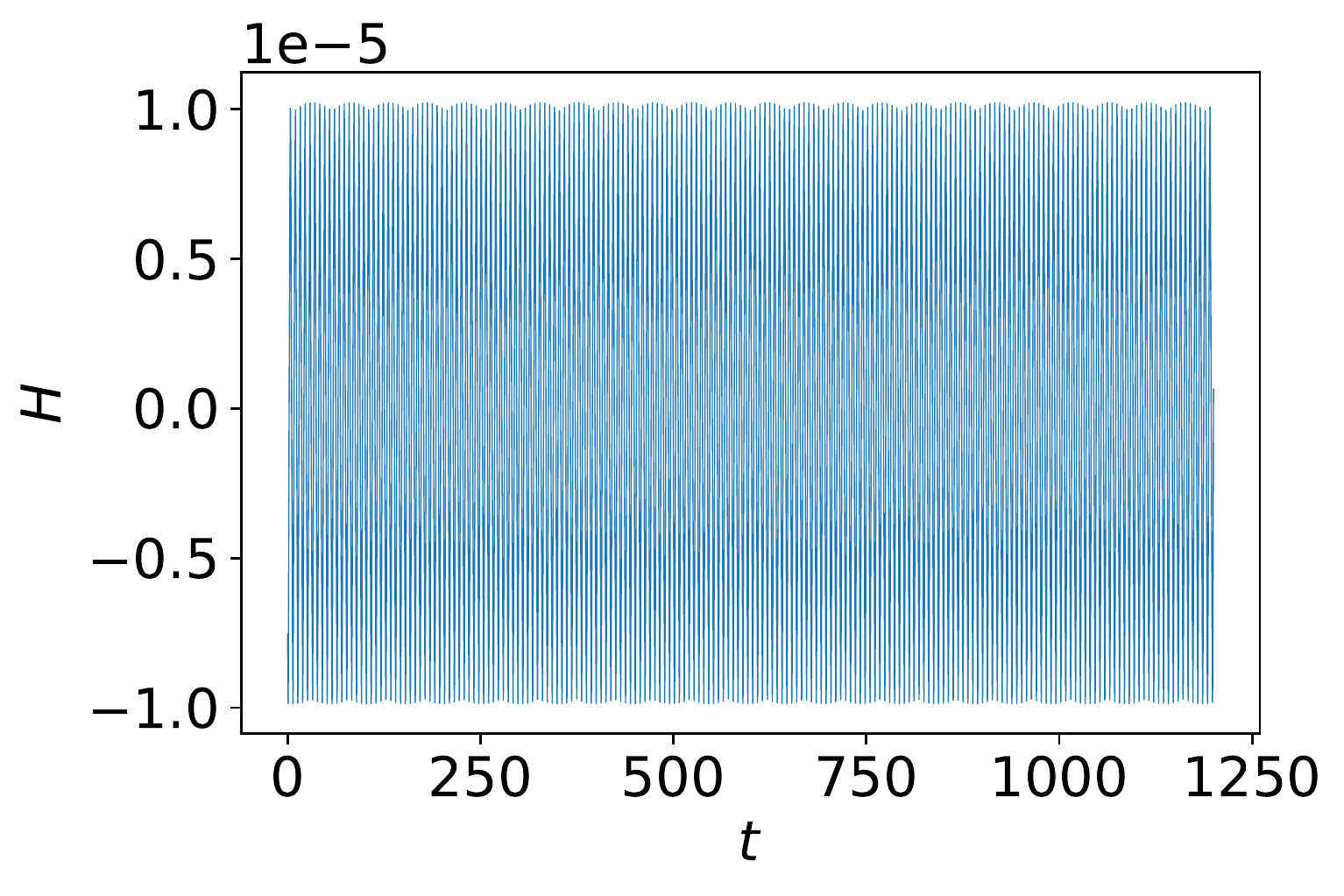}
	}
	\subcaptionbox{As in (c) with more training data}{
		\includegraphics[width=0.4\linewidth]{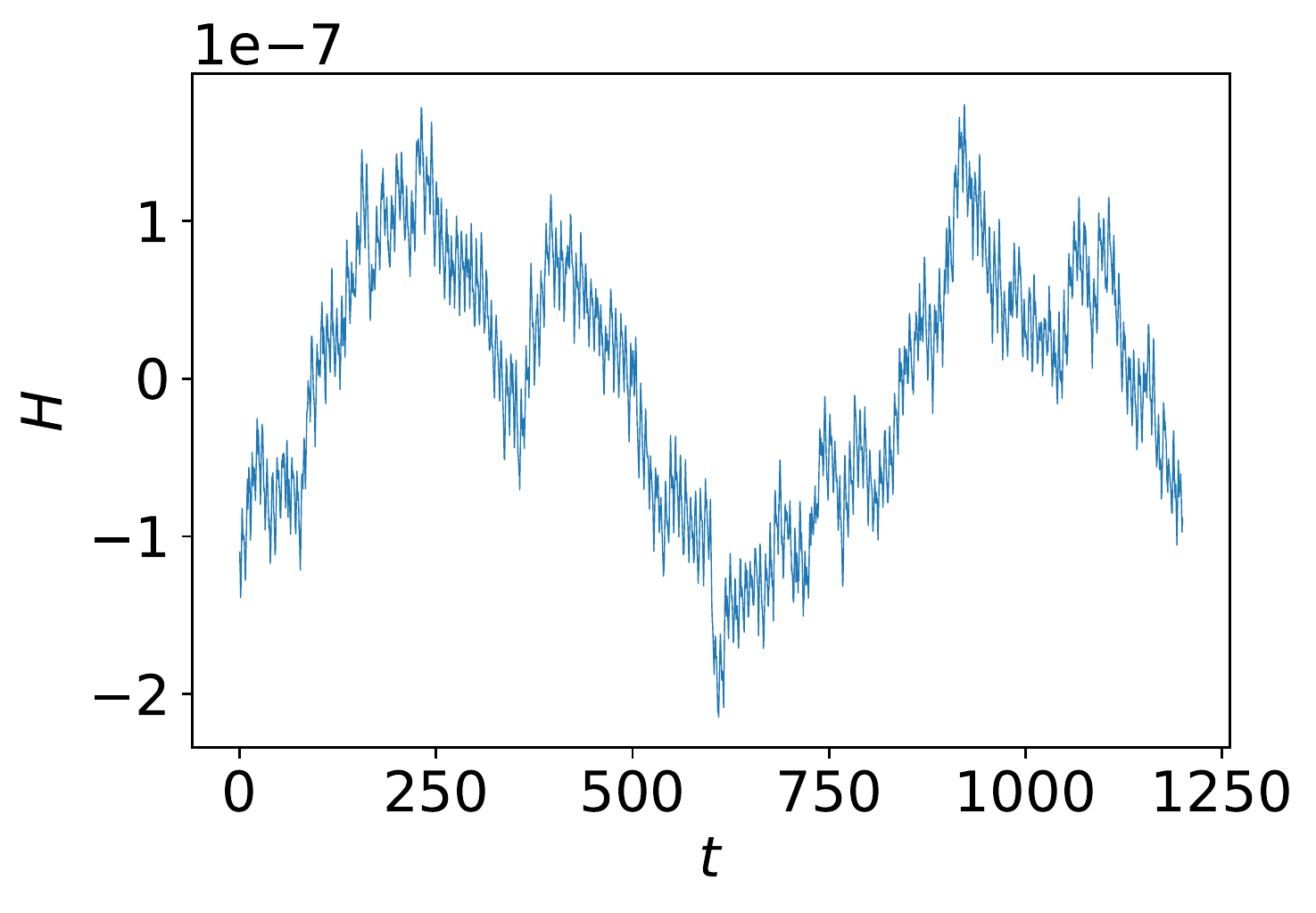}
	}

	\caption{Pendulum experiment with SE, step size $h=0.3$, initialisation at $z_0=(0.4,0)$. The SSI approach is shown in blue, a direct application of SE in grey.
		Figure (a) shows a phase plot of both trajectories as well as a reference solution (dashed) which is exactly covered by the SSI trajectory.
		The Hamiltonian $H$ is plotted along the SSI and SE trajectory in (b). The amplitude of the blue oscillation is significantly smaller than of the grey oscillation.
		(c) depicts the energy of the SSI motion, where the arithmetic mean of $H$ along the trajectory was subtracted.
		When the experiment is repeated with an abundance of training data ($N=700$) the energy error looks like a random walk (d).}\label{fig:PendulumTrj}
\end{figure}

\Cref{fig:PendulumTrj} compares the behaviour of a trajectory obtained using the Symplectic Euler method (SE), when SE is used with SSI and when SE is applied to the exact system $\dot{z}=J^{-1}\nabla H(z)$ directly, which corresponds to strategy 2 in the infinite data limit.
We see that the SSI prediction conserves the exact energy $H$ much better than the strategy 2 prediction. Indeed, the SSI trajectory is visually indistinguishable from a reference solution.
Moreover, if the experiment is repeated with an increased dataset of $N=700$ training points obtained from a Halton sequence, the amplitude of the energy error oscillation further decreases and looks like a random walk which does not leave a band of width $\approx 4 \cdot 10^{-7}$ during the simulation time.
For strategy 2 to have such a small energy error on the given time-interval, we need to decrease the step-sizes by a factor of 10000 since SE is just of first order.
This demonstrates that SSI successfully compensates the discretisation error introduced by SE and that rather large step size $h=0.3$ can be used while maintaining high accuracy.

\begin{figure}
	\subcaptionbox{Trajectories}{
		\includegraphics[height=0.3\linewidth]{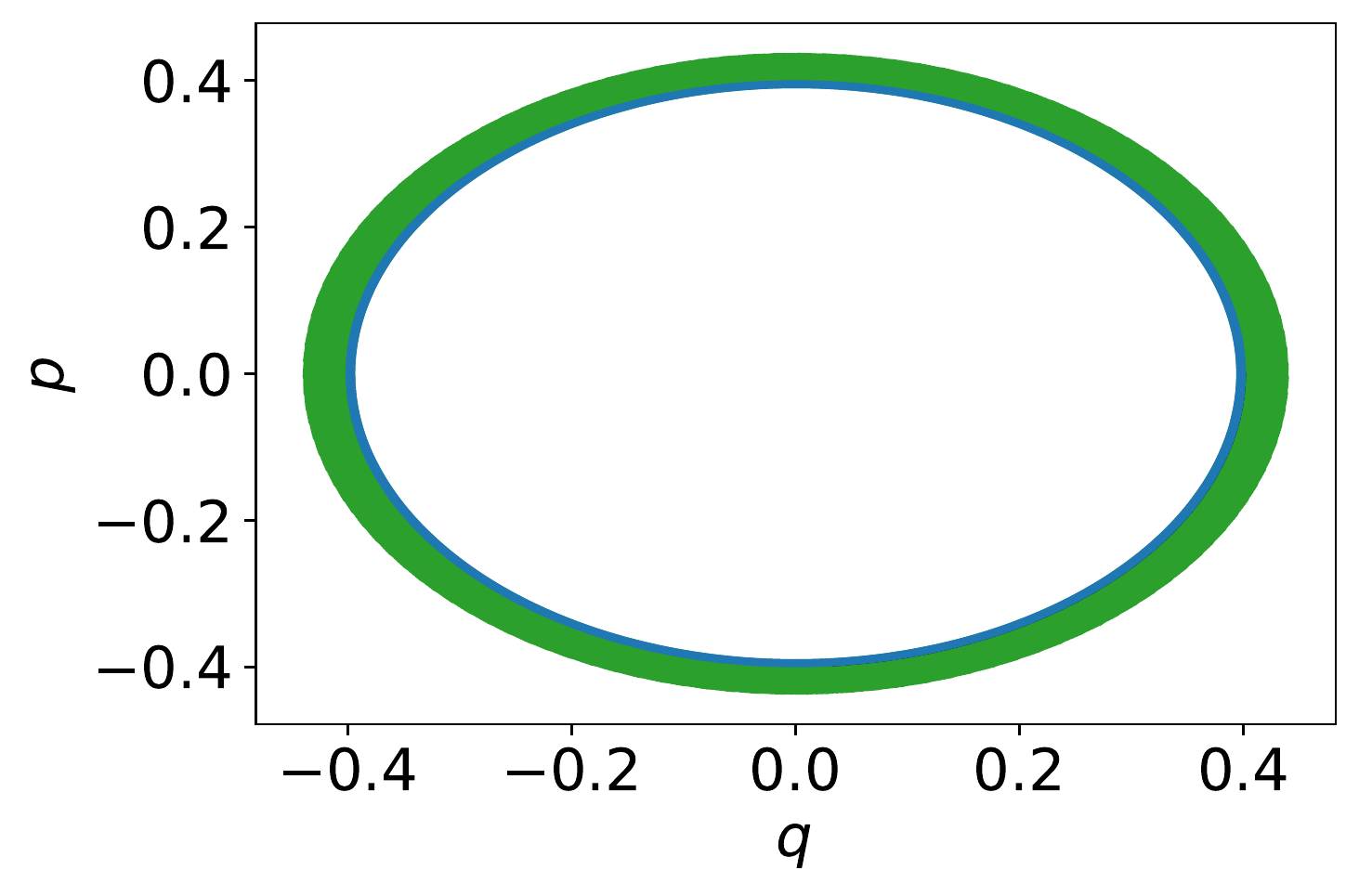}
	}
	\subcaptionbox{Energy}{
		\includegraphics[height=0.3\linewidth]{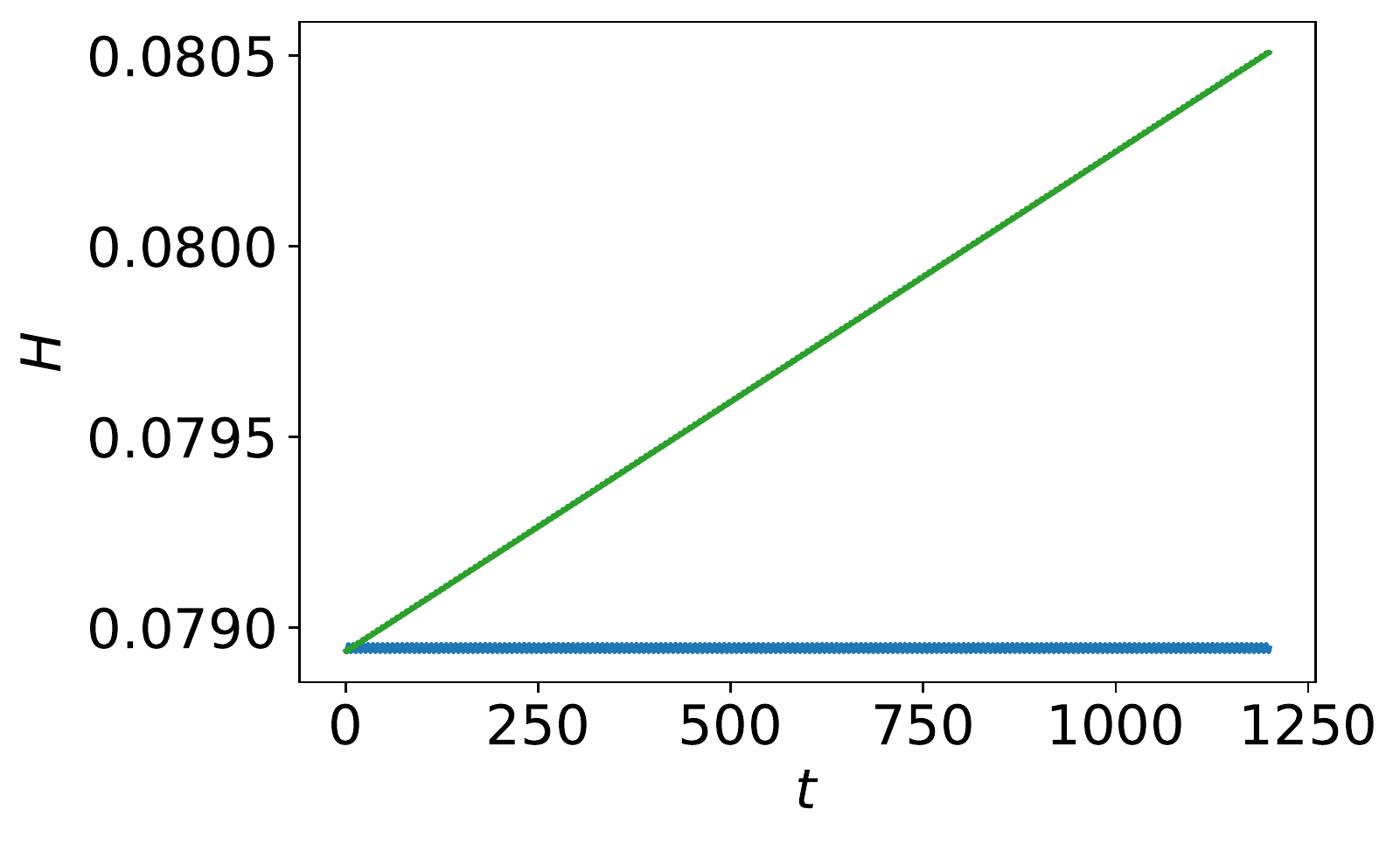}
	}

	\caption{Pendulum experiment. The energy conservation properties of SSI are compared to a GP fitted to the flow map data without incorporation of the Hamiltonian structure. Its trajectory initialised at $z_0=(0.4,0)$ (green) spirals outwards and its energy erroneously grows approximately linearly, while the trajectory from SSI (blue) is periodic and has excellent energy behaviour even on long time intervals.}\label{fig:PendulumSciKit}
\end{figure}

In the following, we continue with the smaller data set $N=160$ and compare to strategy 1, i.e.\ learning the Hamiltonian flow map directly by fitting a Gaussian process to the training data $(Y,\bar{Y})$. \Cref{fig:PendulumSciKit} shows that the learned flow is not energy conserving, which is not surprising as we have not explicitly incorporated Hamiltonian structure. The trajectory starting from $z_0=(0.4,0)$ has a steady energy growth and does not capture the periodicity. This demonstrates that structure preservation is crucial in this example.

\begin{figure}
	\subcaptionbox{Phase portrait}{
		\includegraphics[width=0.44\linewidth]{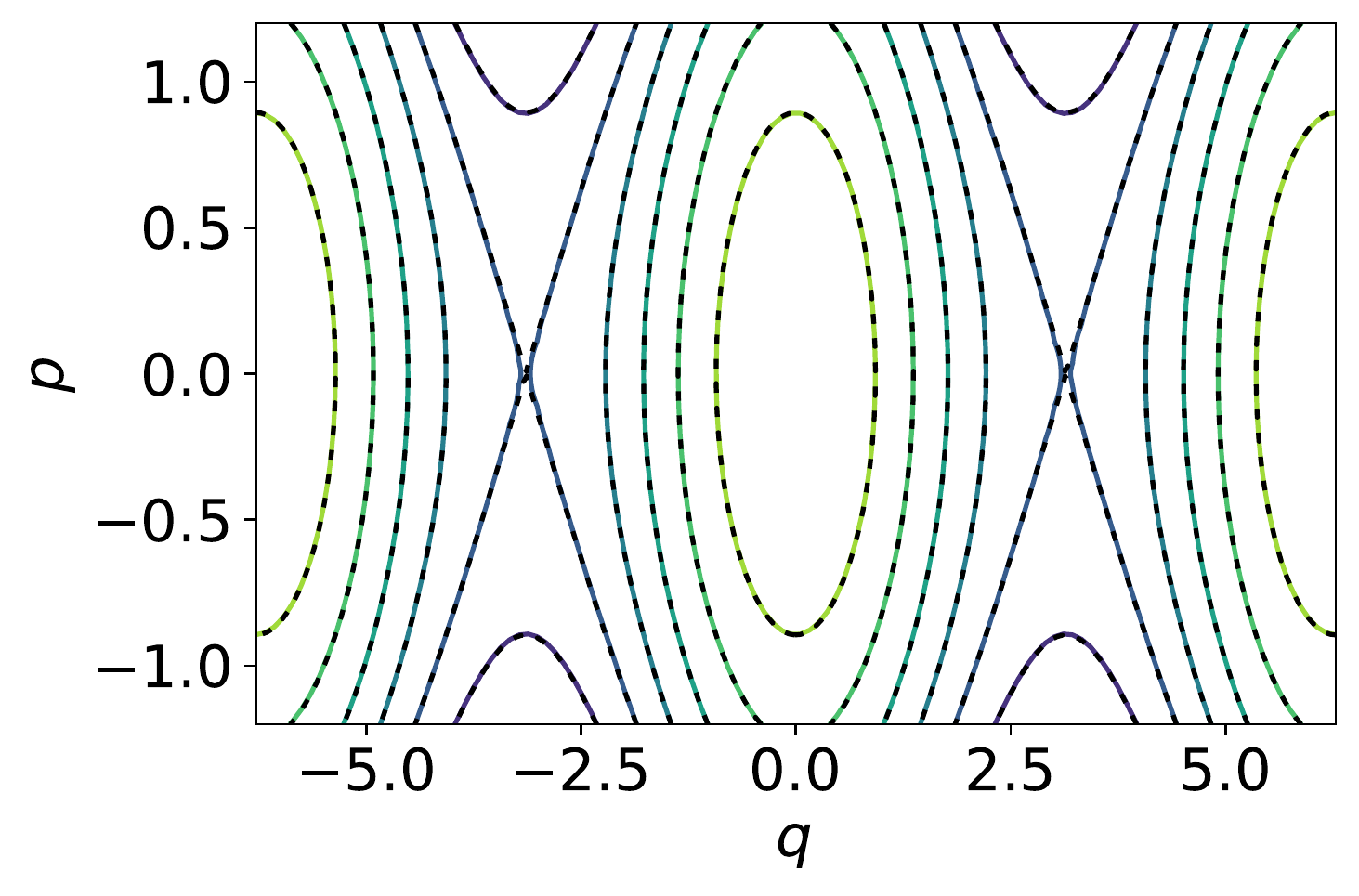}
	}\\
	\subcaptionbox{truncations of $\widetilde{\overline{H}}$}{
		\includegraphics[width=0.44\linewidth]{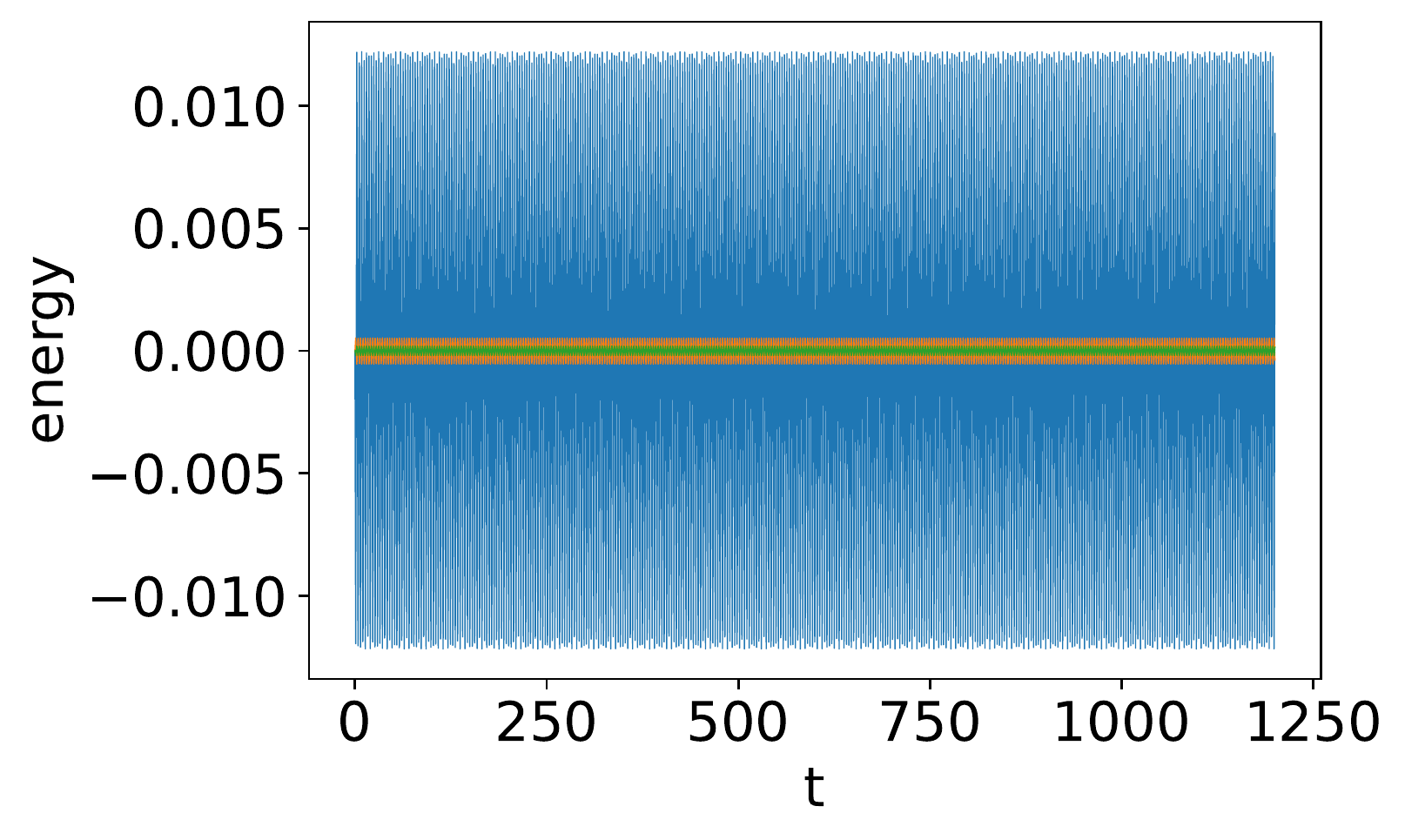}
	}
	\subcaptionbox{truncations of $\widetilde{\overline{H}}$}{
		\includegraphics[width=0.44\linewidth]{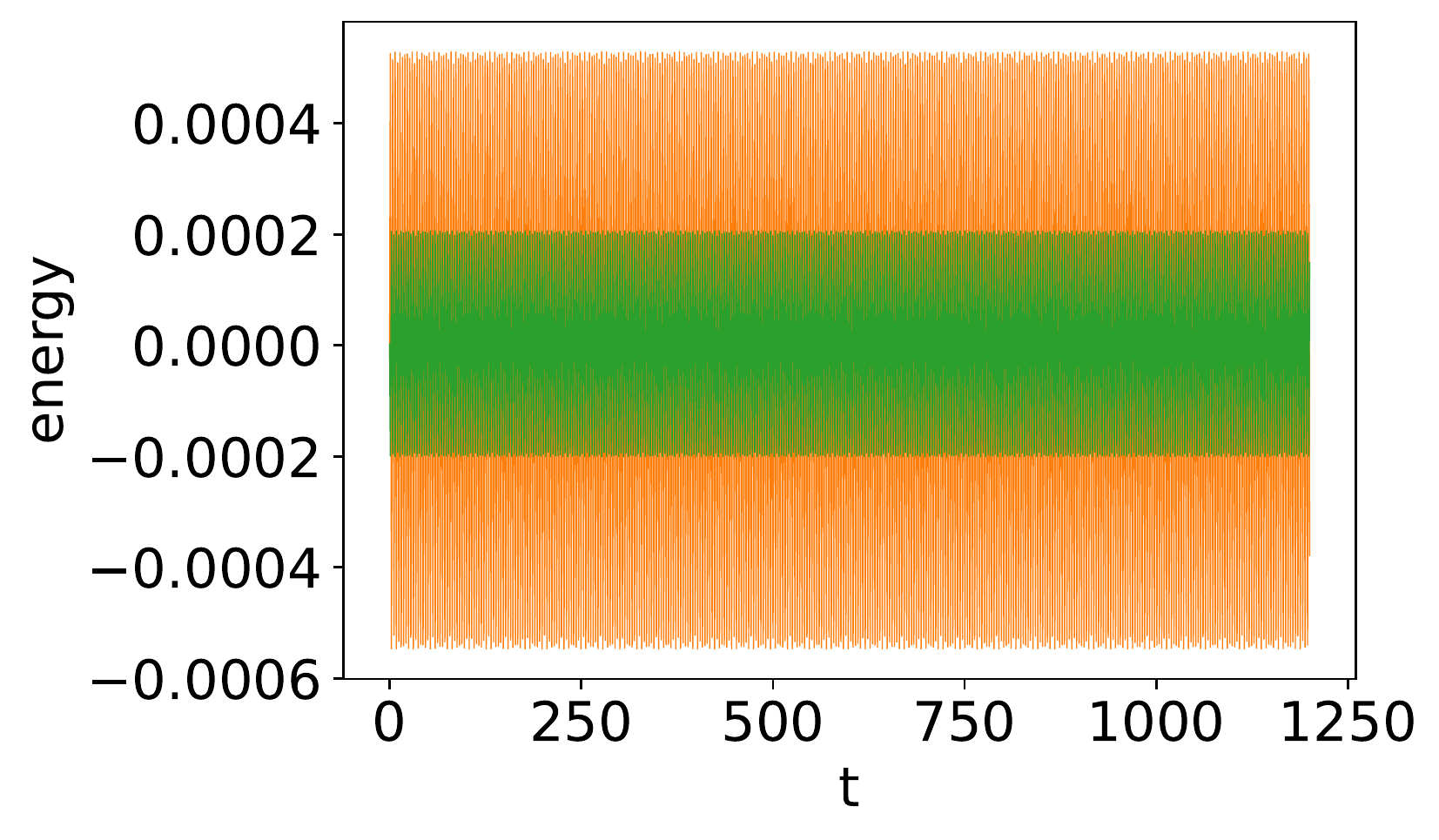}
	}
	
	\caption{Pendulum experiment (SE).
		(a) Level sets of $\widetilde{\overline{H}}^{[2]}$ almost coincides with level sets of $H$ (dashed).
		(b) and (c) show that the truncations of $\widetilde{\overline{H}}$ to 0th (blue), 1st (orange), and 2nd (green) order evaluated along the blue SSI trajectory of \cref{fig:PendulumTrj} (a) are increasingly better preserved. This justifies that a Hamiltonian system for the truncation $\widetilde{\overline{H}}^{[2]}$ describes the numerical motions of the SSI.
	}\label{fig:PendulumMod}
\end{figure}

To analyse qualitative aspects of the SSI scheme, we compute the systems whose exact dynamics coincides with SSI predictions. For this, we apply the backward error analysis formula \eqref{eq:tildeH} to the learned inverse modified Hamiltonian $\overline{H}$ and obtain a power series $\widetilde{\overline{H}}$. 
Several truncations are evaluated along the trajectory in \cref{fig:PendulumMod}. As the oscillations of the second truncation $\widetilde{\overline{H}}^{[2]}$ occur within a band of width $4 \cdot 10^{-4}$, a truncation to 2nd order describes the Hamiltonian that governs the numerical dynamics sufficiently well for our purposes.
Indeed, a contour plot of $\widetilde{\overline{H}}^{[2]}$ on $M$ is visually indistinguishable from a contour plot of the exact $H$. Here we have used a uniform mesh with 120x120 points.
A Hamiltonian is only defined up to a constant by the system's motion. Rather than measuring the $L^{2}$ distance of $H$ and $\widetilde{\overline{H}}^{[2]}$ over $M$, we compute the standard deviation of $H_{\mathrm{diff}} = H - \widetilde{\overline{H}}^{[2]}$ for a uniform distribution on the training domain $M$, i.e.\
\[
\sigma(H_{\mathrm{diff}}) 
= \sqrt{\mathbb{E}[(H_{\mathrm{diff}}-\mathbb{E}[H_{\mathrm{diff}}])^2]}, \, \text{with} \,
\mathbb{E}[f] = \frac{\int_M f \d \nu}{\d \nu(M)}
\]
and with $\d \nu$ denoting the Lebesgue measure on $M$. Approximating $M$ using a 120x120 uniform mesh, we obtain $\sigma(H_{\mathrm{diff}}) \le 5.2 \cdot 10^{-4}$. This explains the excellent energy conservation of SSI motions and shows that qualitative aspects such as periodic motions and even the behaviour close to the separatrix are guaranteed, since the phase portrait of the system governing the SSI motions is close to the exact phase portrait. This is in particular important for long-term simulations.
Furthermore, this demonstrates that SSI can be used for system identification.


\begin{figure}
	\subcaptionbox{MP}{
		\includegraphics[width=0.44\linewidth]{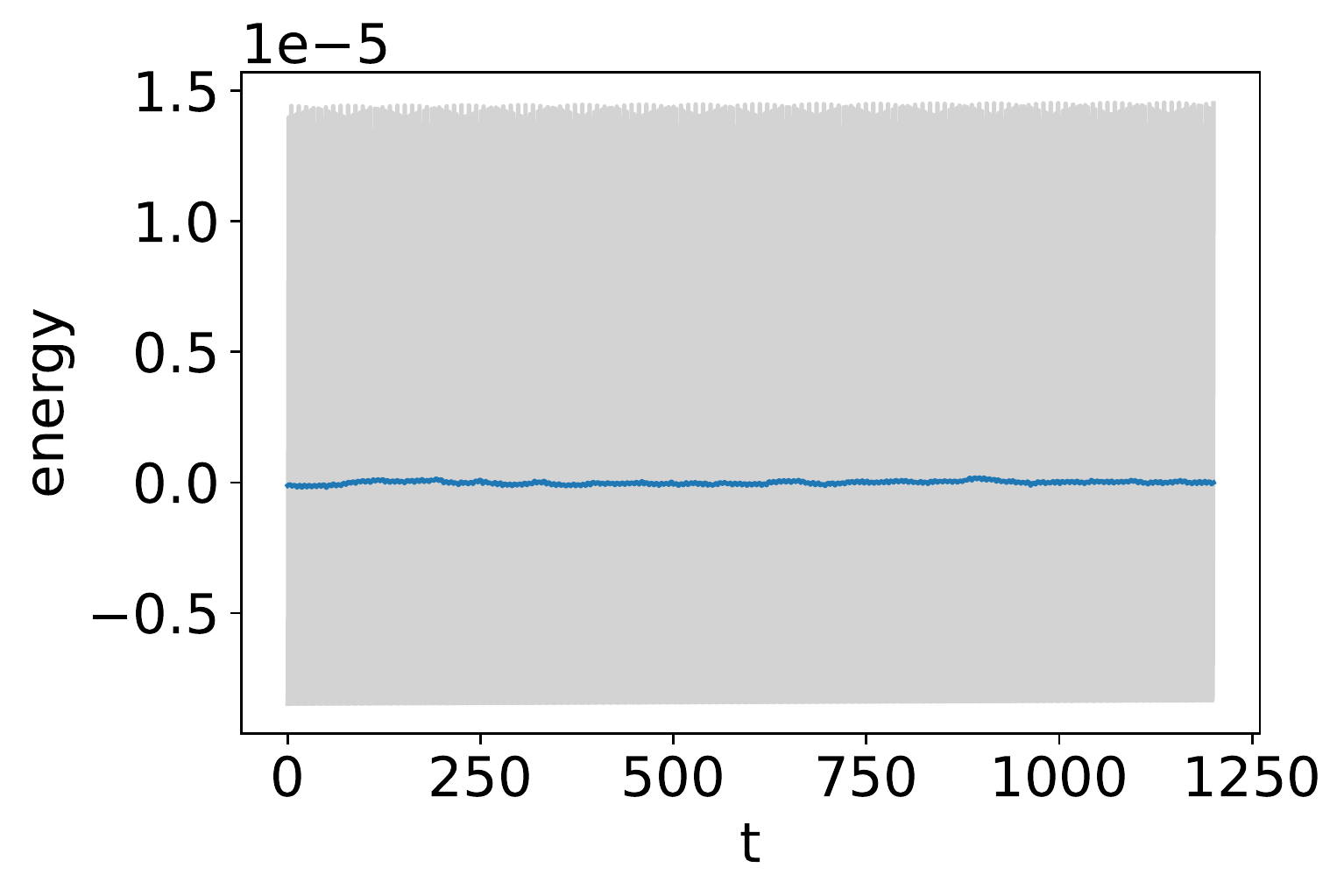}
	}
	\subcaptionbox{SSI}{
		\includegraphics[width=0.44\linewidth]{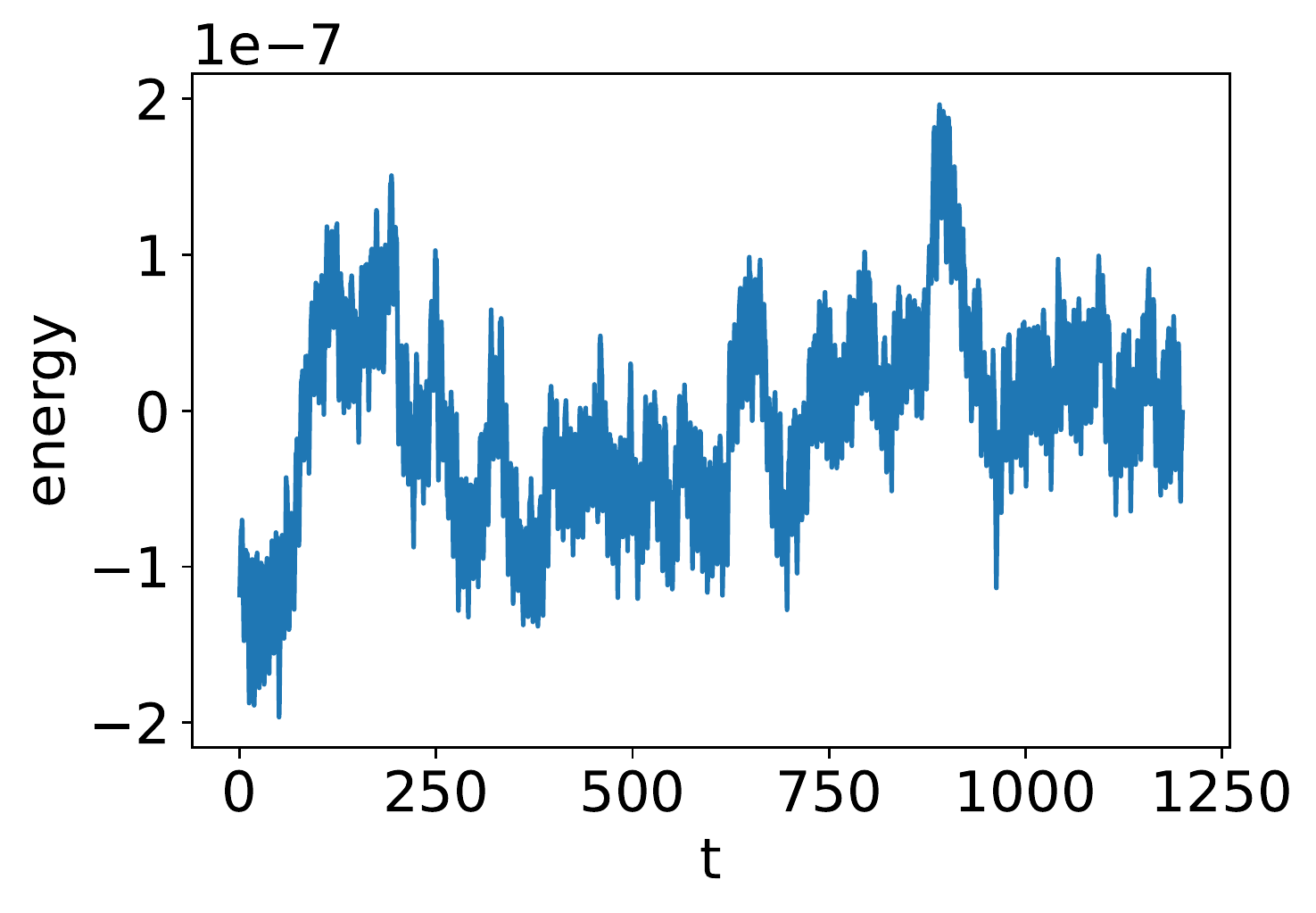}
	}
	\caption{Pendulum experiment (MP).
		The energy behaviour of MP (gray) applied to the exact system is oscillatory while the SSI approach shows energy conservation up to a random walk (blue). In all plots, the mean value of the energy trajectory has been subtracted before plotting.
	}\label{fig:PendulumTrjMD}
\end{figure}

The experiments are repeated with the second order accurate Implicit Midpoint rule (MP) using the same parameters but $N=400$ data points in the training process.
Phase plots of the trajectory computed with SSI and with a direct application of MP to the exact system $\dot{z} = J^{-1} \nabla H(z)$ (strategy 2 in the infinite data limit) both visually coincide with an exact trajectory (not shown).
However, while the energy error of the strategy 2 approach is oscillatory with amplitudes bounded by $3 \cdot 10^{-5}$, the energy error of the SSI solution shows a random walk within a band of width $4 \cdot 10^{-7}$ on long time scales (\cref{fig:PendulumTrjMD}). For a similar energy error on the given time-interval we would need to decrease the step size in the strategy 2 approach by a factor of approximately 9. Moreover, the error shows a linear trend such that the step size needs to be decreased further when simulation time is increased (not shown).

\begin{figure}
	\subcaptionbox{0th/1st (orange) and 2nd/3rd (green) truncation of $\widetilde{\overline{H}}$}{
		\includegraphics[width=0.44\linewidth]{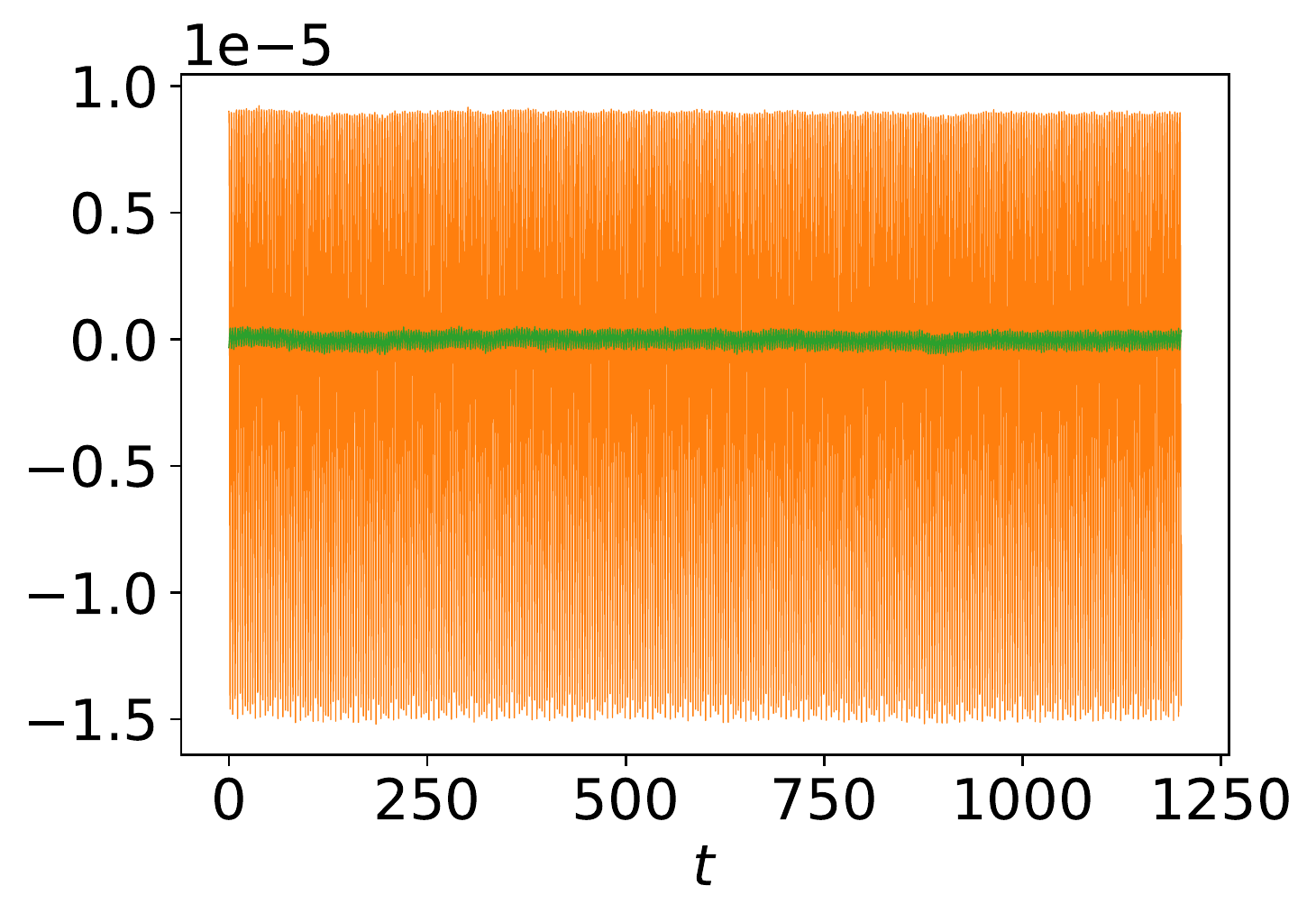}
	}
	\subcaptionbox{2nd/3rd truncation of $\widetilde{\overline{H}}$}{
		\includegraphics[width=0.44\linewidth]{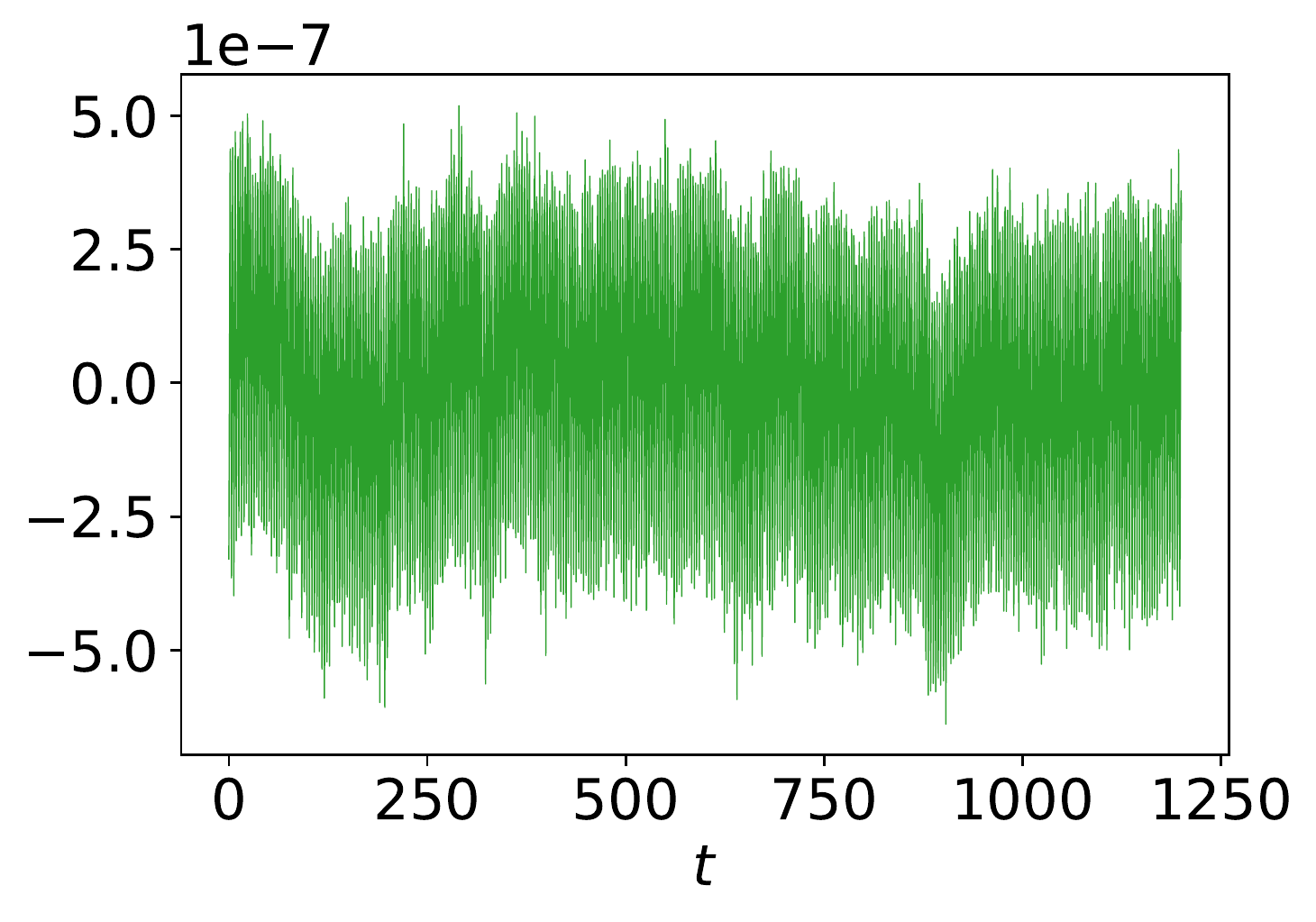}
	}
	\caption{Pendulum experiment (MP). Repetition of the experiment of \cref{fig:PendulumMod} with MP instead of SE. Notice that the 0th and 1st as well as the 2nd and 3rd truncation of $\widetilde{\overline{H}}$ coincide due to the symmetry of MP.
	}\label{fig:PendulumModMD}
\end{figure}

Similarly as in the experiments with SE, the truncation $\widetilde{\overline{H}}^{[2]}$ computed with the backward error analysis formula \eqref{eq:MDinvmodH} describes the numerical dynamics accurately as it is up to small errors a conserved quantity of the numerical flow (\cref{fig:PendulumModMD}). Notice that the 0th and 1st order truncation and the 2nd and 3rd order truncations of $\widetilde{\overline{H}}$ coincide since the midpoint rule is symmetric. The standard deviation of $H_{\mathrm{diff}}=H-\widetilde{\overline{H}}^{[2]}$ (computed as before) is smaller than $\sigma(H_{\mathrm{diff}}) \le 9.4 \cdot 10^{-4}$. Again, this shows that SSI replicates qualitative aspects of the exact system even in long term simulations.

We see that using the Symplectic Euler scheme or the Implicit midpoint rule in the SSI technique yields similar behaviour and there appears to be no gain from using a higher order method for SSI.


\subsection{H\'enon--Heiles system}

\begin{figure}
	\includegraphics[width=0.44\linewidth]{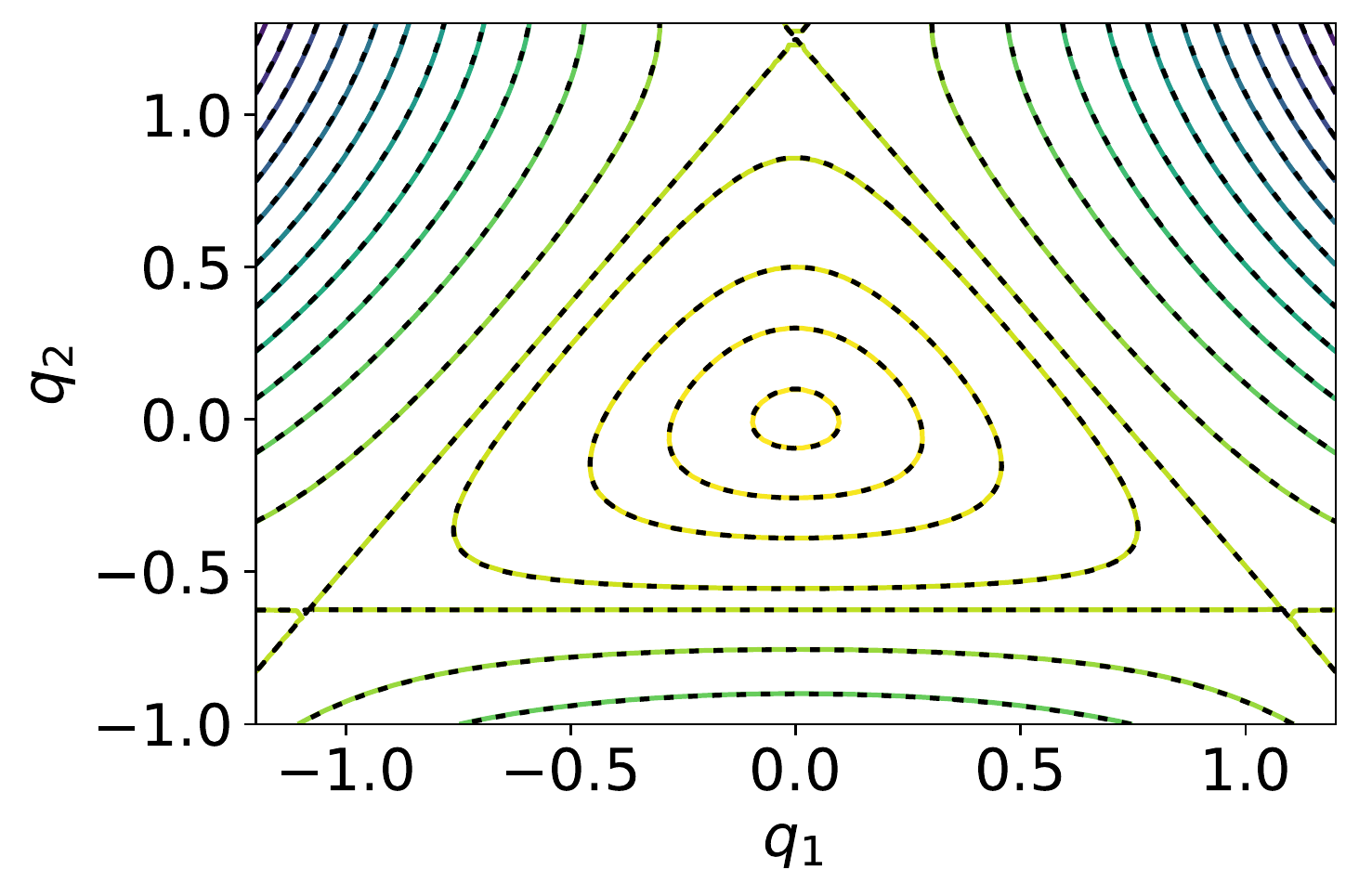}
	\caption{Contour plot of the exact H\'enon--Heiles potential $V$ (dashed) and the potential $\tilde{\bar{V}}^{[2]}(q) := \widetilde{\overline{H}}^{[2]}(q,(0,0))$ computed from the learned inverse modified Hamiltonian $\overline{H}$. The level sets of $\tilde{\bar{V}}^{[2]}$ and $V$ are not distinguishable in this plot.}\label{fig:HenonHeils}
\end{figure}

The H\'enon--Heiles system is a Hamiltonian system with
\[
H(q,p) = \frac{1}{2}\|p\|^2 + V(q), \quad V(q)= \frac{1}{2} \|q\|^2 + \mu \left(q_1^2 q_2-\frac{q_2^3}{3}\right).
\]
We set the parameter to $\mu=0.8$. A contour plot of $V$ is shown in \cref{fig:HenonHeils}.
The system has bounded as well as unbounded motions. More precisely, bounded connected components of level sets of $V$ and $H$ correspond to values of $V$ or $H$ within the interval $I=[0,\frac{1}{6 \mu^2})$. All connected components of level sets to values in $\R \setminus I$ are unbounded.
Accurate preservation of energy is crucial, when integrating trajectories, since trajectories on close-by energy level sets can show very different long-term behaviour. Classically, this forces one to use very small time steps in numerical integration schemes. However, as SSI compensates discretisation errors, we can use a moderate time-step.
%
%
We repeat the experiments from \cref{sec:Pendulum} with $h=0.1$, $N=800$, $n_{\mathrm{LP}}=600$, kernel parameter $e = 5$, and the four dimensional tesseract $M=[-1,1]^4$ as domain for training data.

\begin{figure}
		\includegraphics[width=0.44\linewidth]{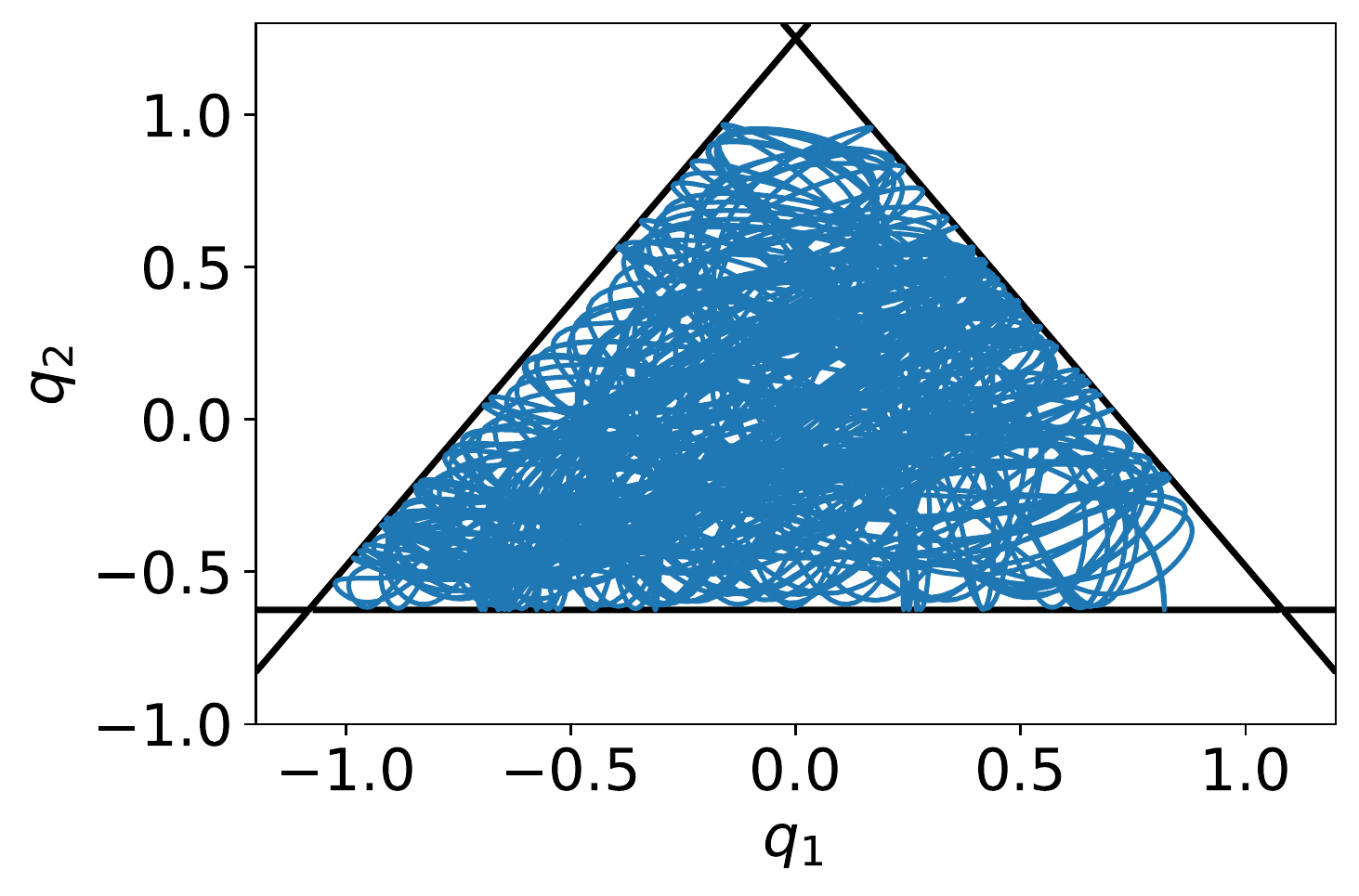}
		\includegraphics[width=0.44\linewidth]{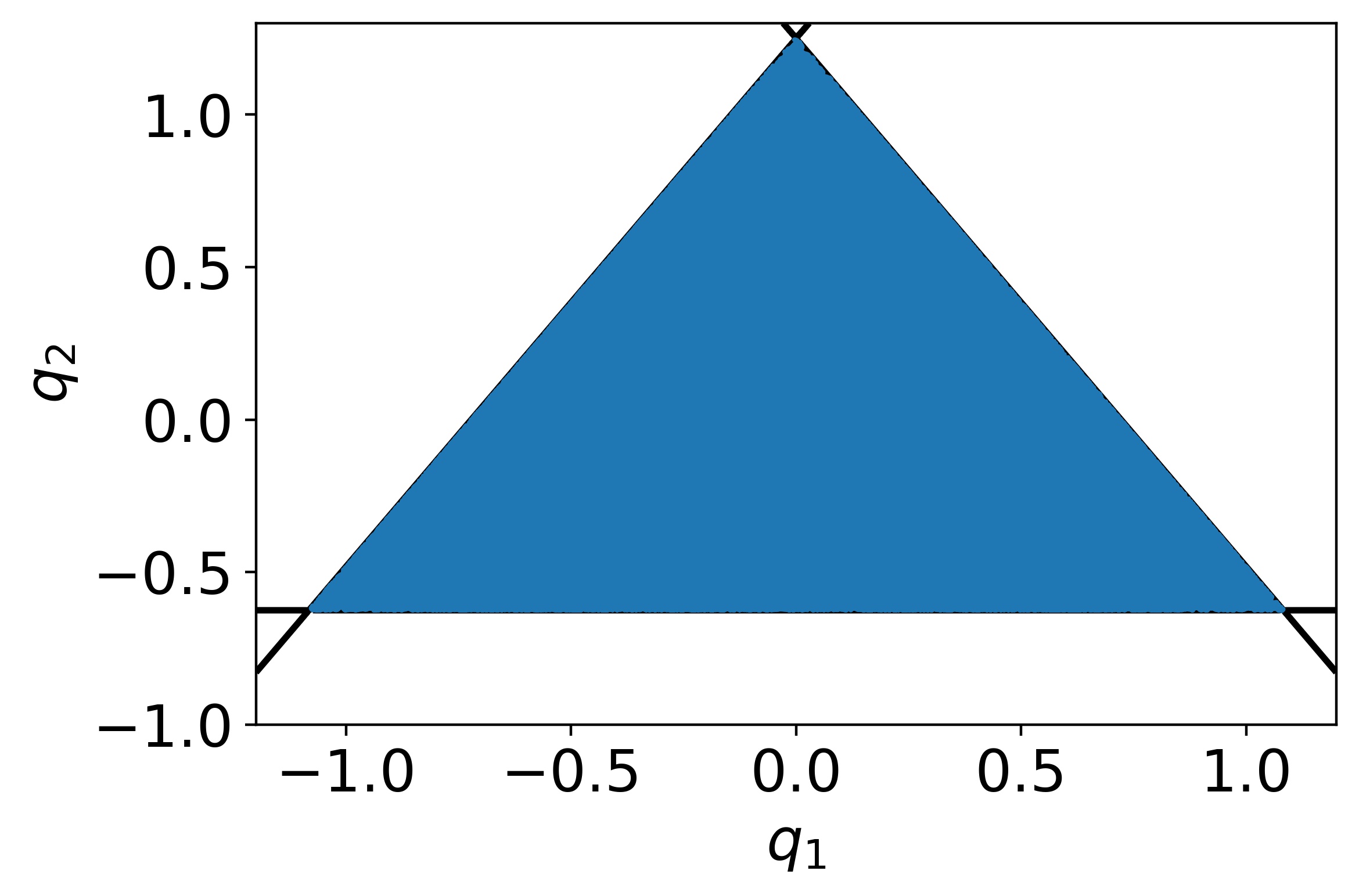}
	\caption{H\'enon--Heiles experiment. Phase plot of trajectory from SSI initialised at $z_0$ computed up to time $T=800$ and $T=40000$. The motion (blue) densely fills the area bounded by the level set $V^{-1}(q_0)$ (black). }\label{fig:HenonHeilsTrj}
\end{figure}

\begin{figure}
	\begin{center}
	\subcaptionbox{Trajectory SE}{
		\includegraphics[width=0.44\linewidth]{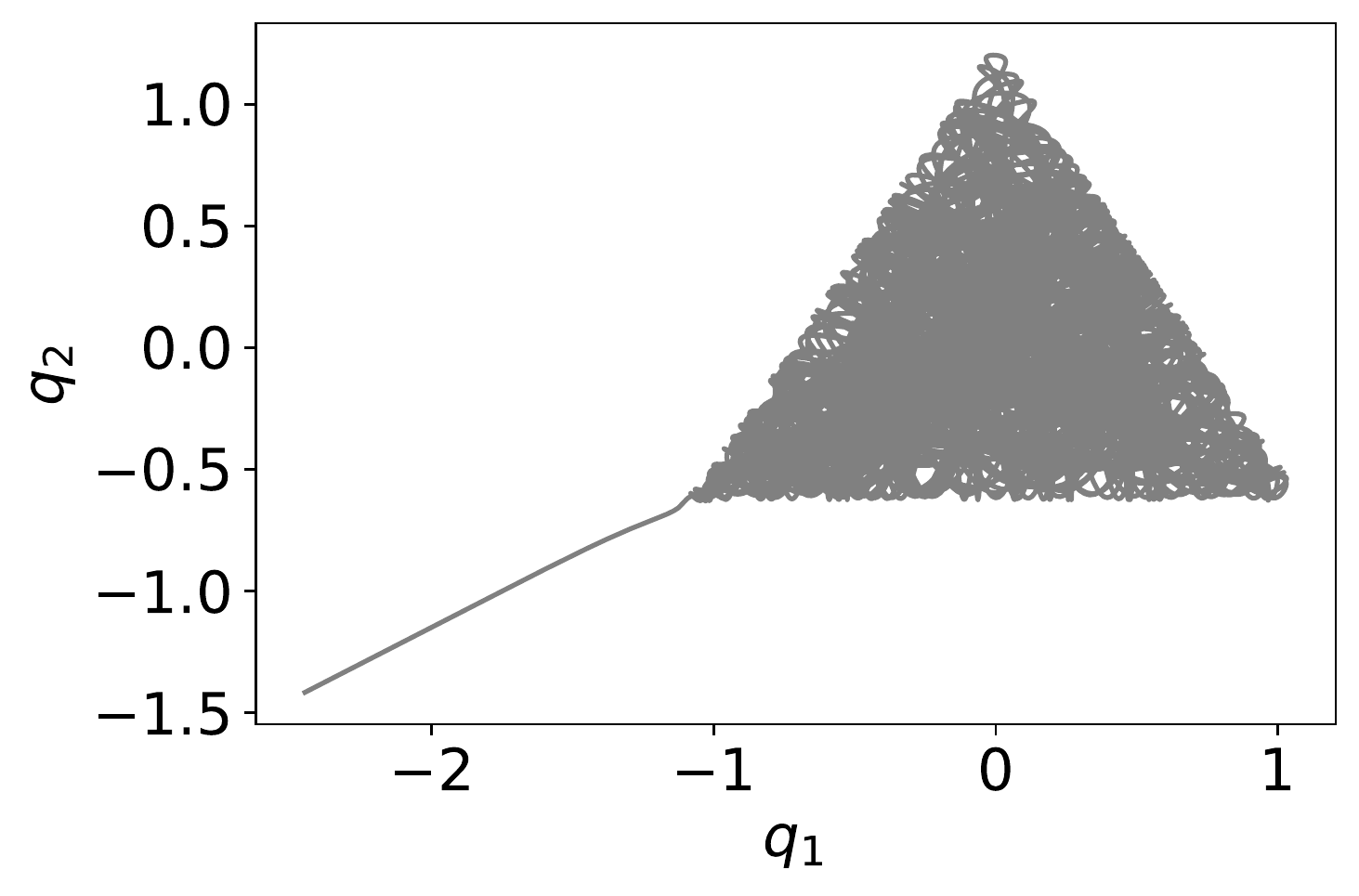}
	}
	\subcaptionbox{Energy SE/SSI}{
		\includegraphics[width=0.44\linewidth]{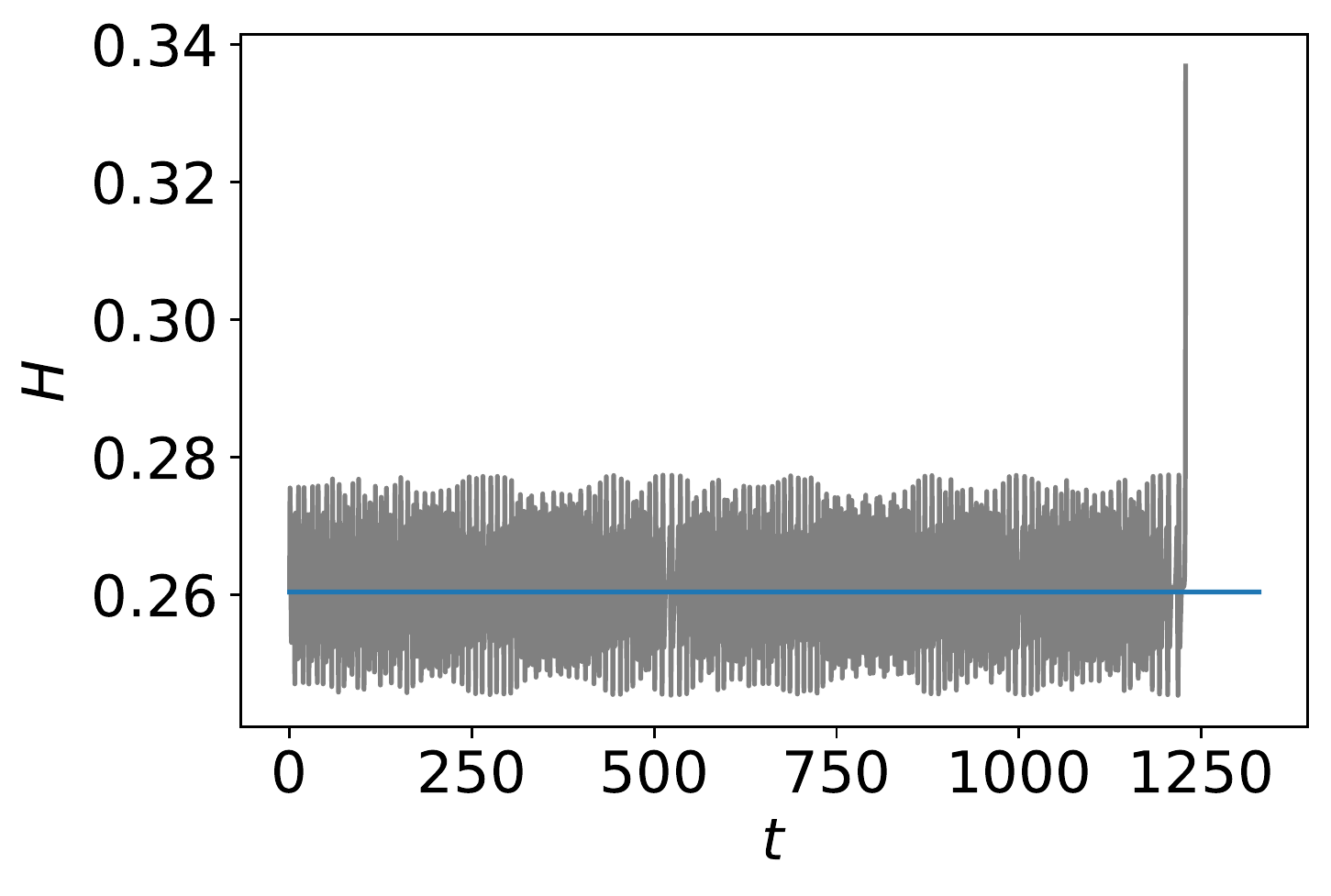}
	}
	\subcaptionbox{Energy SSI}{
		\includegraphics[width=0.44\linewidth]{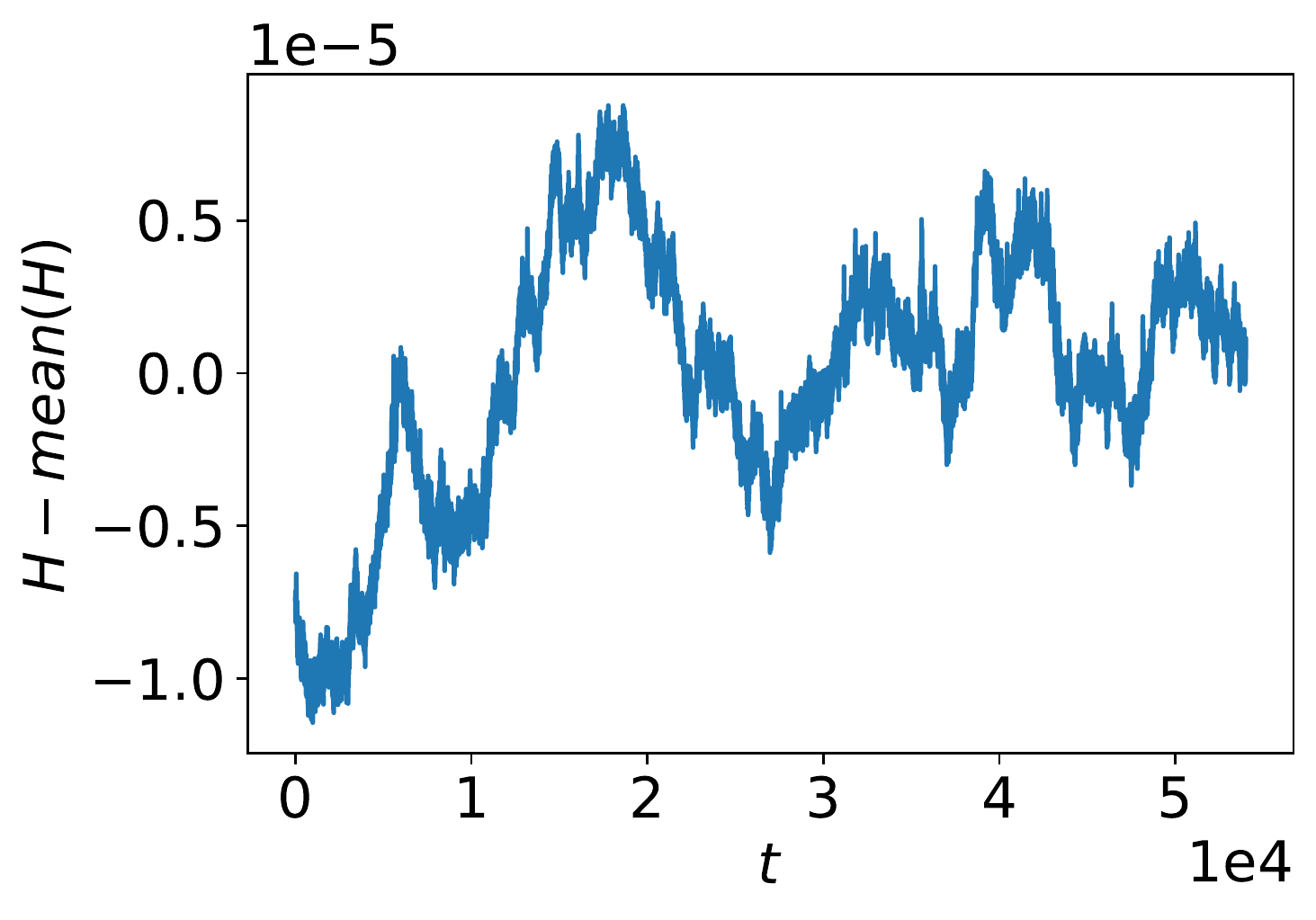}
	}
	\subcaptionbox{SE $h/10000$}{
	\includegraphics[width=0.44\linewidth]{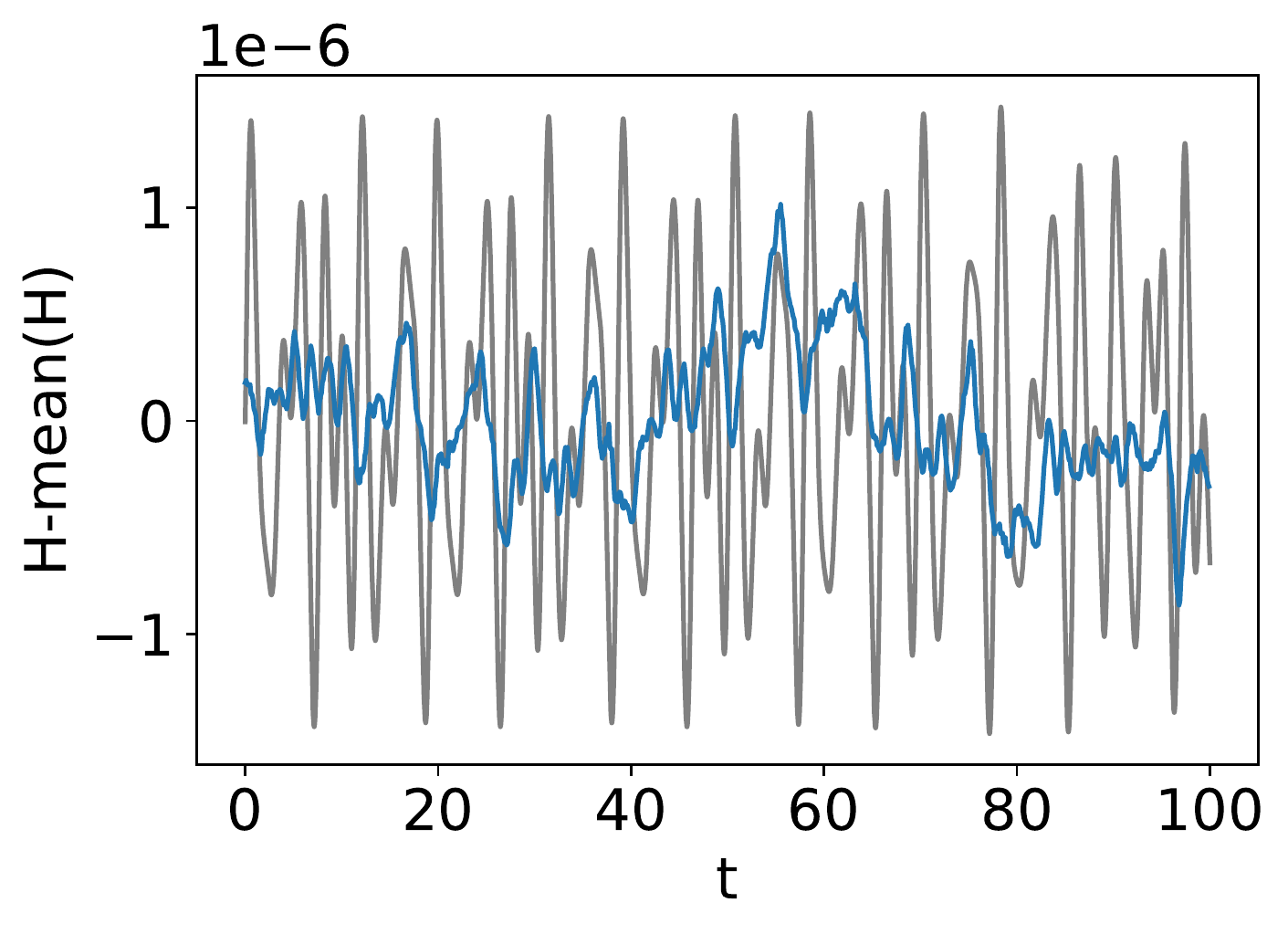}
}

	\end{center}	
	
	\caption{H\'enon--Heiles experiment. The trajectory obtained from an application of SE to the exact differential equation $\dot{z}=J^{-1} \nabla H(z)$ (grey) erroneously leaves its energy level set and diverges. (a) shows the trajectory up to time $t=1229.5$. (b) compares energy conservation of SSI (blue) with SE (grey). (c) shows the random walk like long-term energy behaviour of SSI. (d) compares energy conservation of SSI and SE on a smaller time interval, where SE's time-step was decreased to $h/10000$.
	}\label{fig:HenonHeilsTrjH}
\end{figure}

The exact motion initialised at $z_0=(q_0,p_0)=((0.675499,0.08),(0,0))$ lies on a bounded connected component of its energy level set $H(z_0)$. While the trajectory predicted by SSI based on SE correctly captures this behaviour (see \cref{fig:HenonHeilsTrj}), a direct application of SE to the exact system $\dot{z}=J^{-1} \nabla H(z)$ with the same step size gives a trajectory which is unbounded. The reason is that $H(z_0)$ is very close to the critical value $\frac{1}{6 \mu^2}$. Indeed $|H(z_0)-\frac{1}{6\mu^2}| \approx 6.2 \cdot 10^{-7}$. SE conserves energy only up to an oscillation within a band of width $\approx 0.03$ such that
the numerical trajectory repeatedly jumps to unbounded energy level sets and eventually escapes (\cref{fig:HenonHeilsTrjH}). In contrast, the trajectory obtained using SSI has a significantly better energy behaviour looking like a random walk, which is bounded within a band of width $\approx 2 \cdot 10^{-5}$ until time $T=5  \cdot 10^4$. While this does not exclude that the trajectory escapes from an energy perspective, it makes the escape less likely. Indeed, no escape occurs during the simulation time.
The step size would need to be decreased by a factor of at least 10000 to obtain comparable energy conservation properties with the strategy 2 approach.

This demonstrates that compensating discretisation errors using the SSI framework can be highly beneficial, when accurate energy preservation is important.

\begin{figure}
	\subcaptionbox{Trajectory}{
		\includegraphics[width=0.44\linewidth]{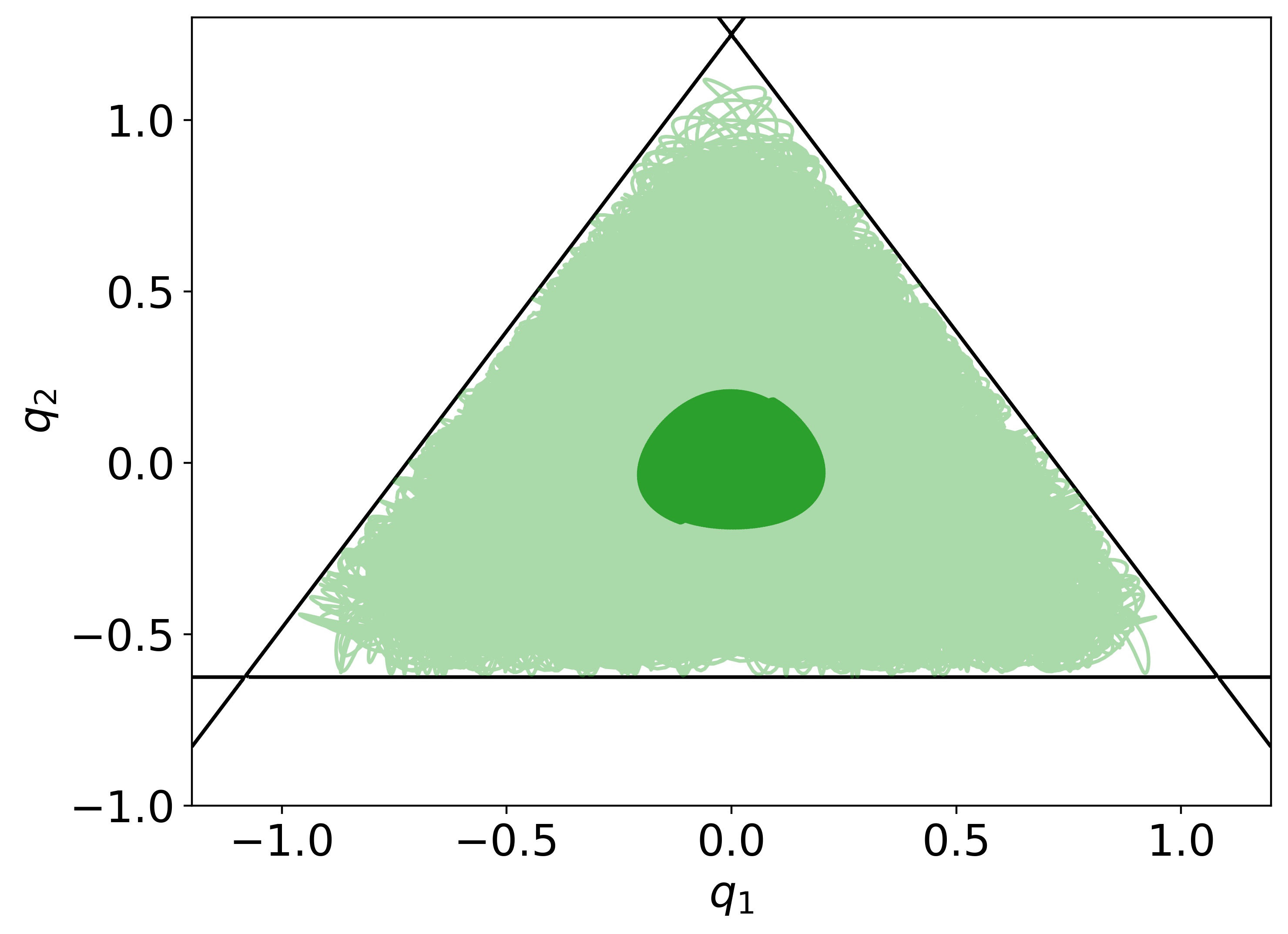}
	}\label{subfig:phaseSciKit}
	\subcaptionbox{Energy}{
		\includegraphics[width=0.44\linewidth]{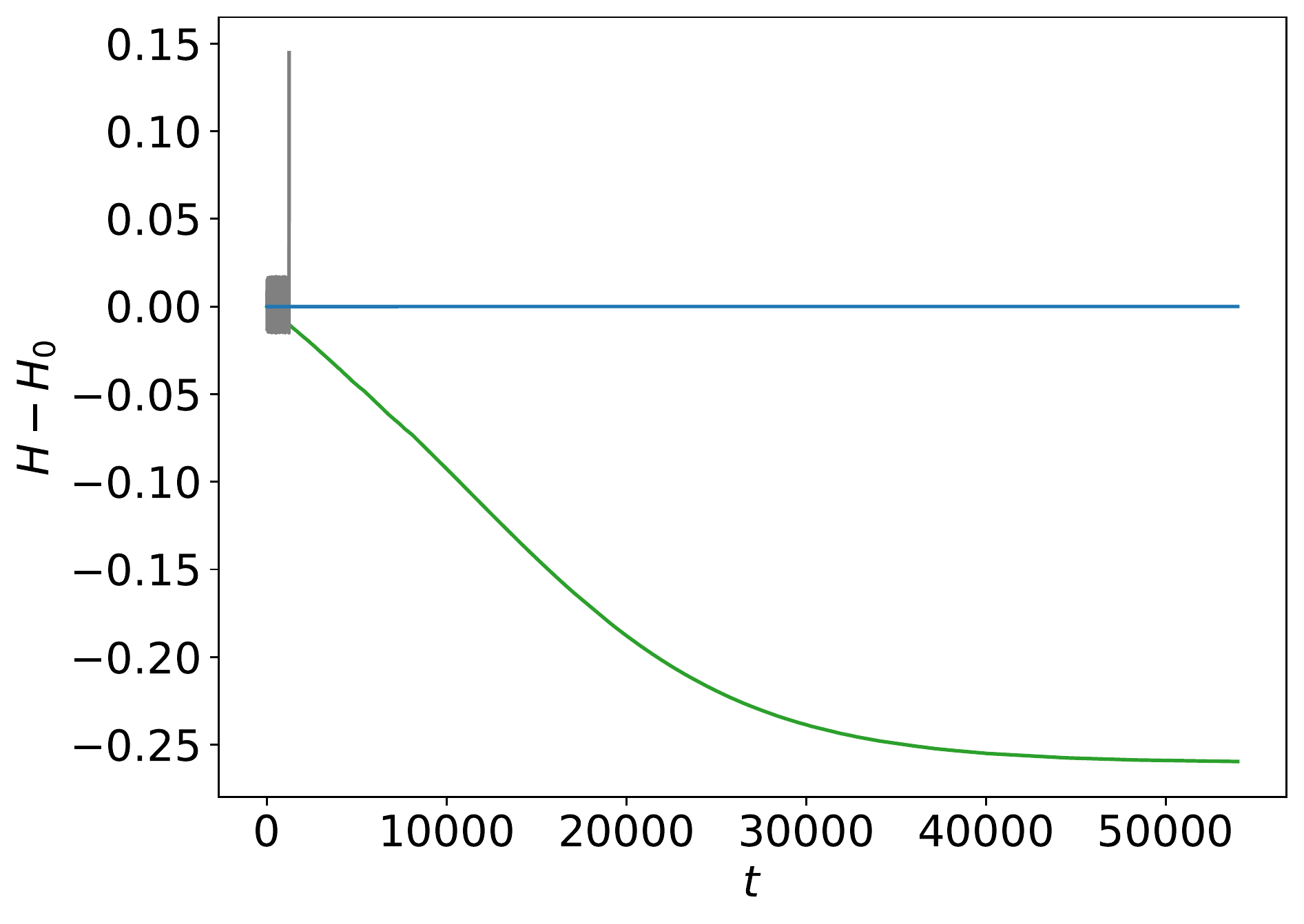}
	}\label{subfig:energySciKit}
	\caption{A motion predicted from a GP fitted to the flow map data $(Y,\bar{Y})$ converges to an artificial attractor and fails to densely fill the area bounded by $V^{-1}(q_0)$. (a) shows a plot of the trajectory in light green up to time $T_1=30000$ and in darker green from time $T_1$ to $T=54000$. (b) compares the energy behaviour of the trajectory from (a) (green), from SSI (blue), and from the direct SE integration (\cref{fig:HenonHeilsTrjH} a) (grey). Only SSI shows the correct long-term energy behaviour.
	}\label{fig:HenonHeilsSciKit}
\end{figure}

For comparison with strategy 1, we learn the flow map of the system directly rather than the inverse modified Hamiltonian $\overline{H}$. This is done by fitting a GP directly to the training data $(Y,\bar{Y})$, again using radial basis functions as kernels \eqref{eq:RBF} and fitting the hyperparameters $k_c, e$ using marginal likelihood estimation.
\Cref{fig:HenonHeilsSciKit} shows a trajectory initialised at $z_0 = (q_0,(0,0))$. It fails to densely fill the area bounded by $V^{-1}(q_0)$. Instead, it monotonously loses energy and artificially converges to the origin. Moreover, access to the learned flow map does not directly reveal information on the Hamiltonian structure of the system, which could be used for system identification.

\begin{figure}
	\subcaptionbox{0th, 1st, 2nd truncation}{
		\includegraphics[height=0.3\linewidth]{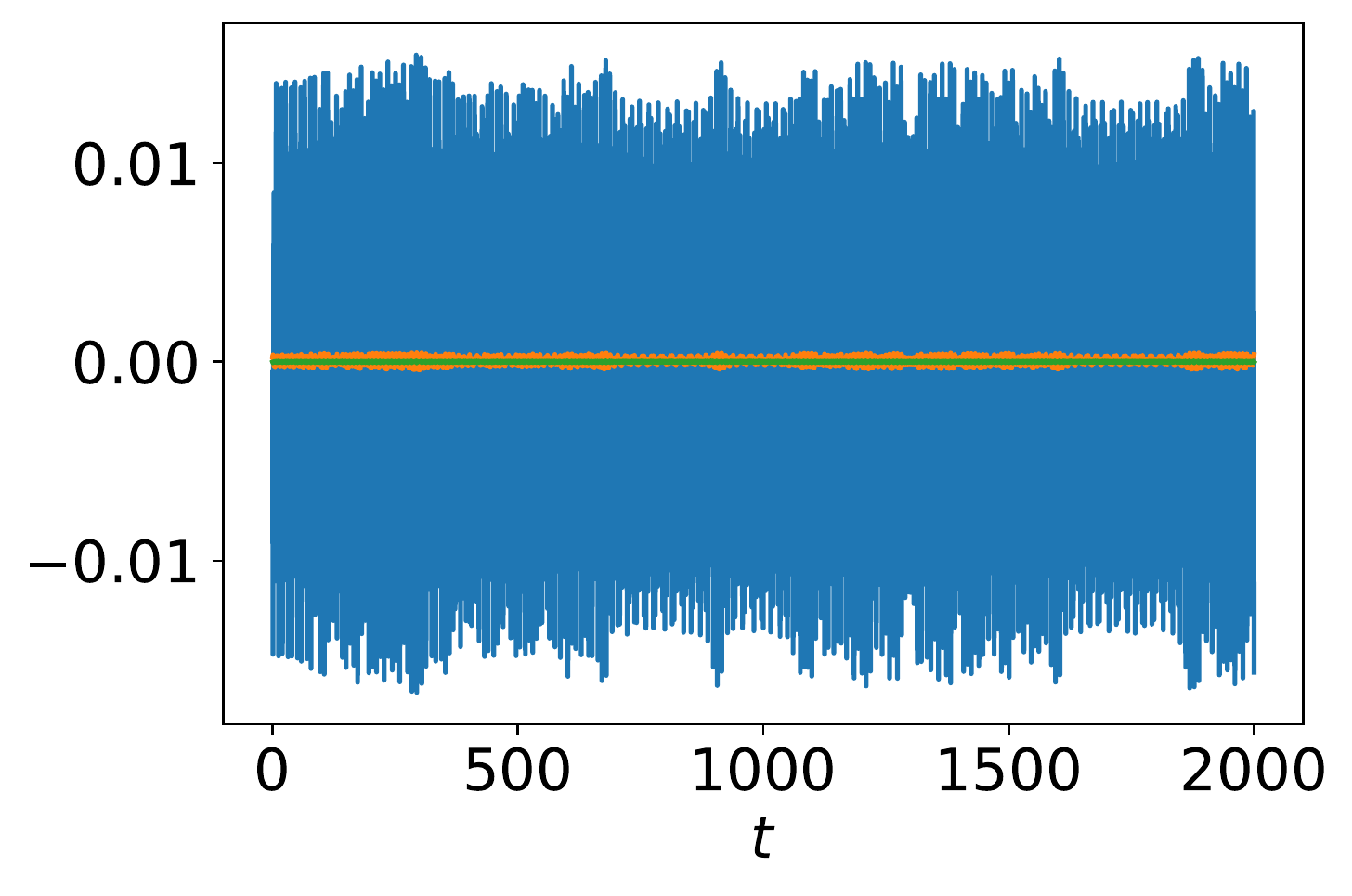}
	}\\
	\subcaptionbox{1st, 2nd truncation}{
		\includegraphics[height=0.3\linewidth]{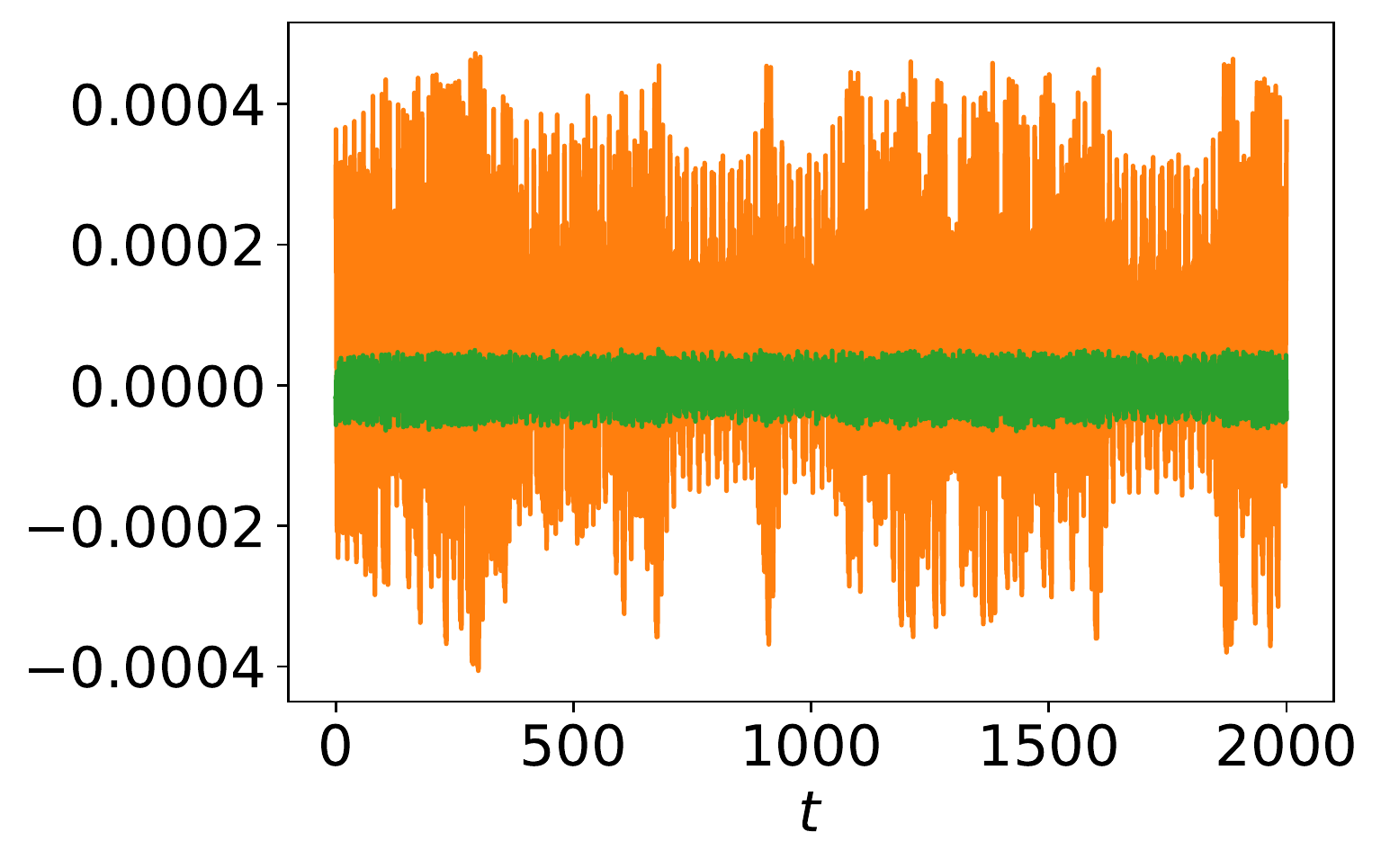}
	}
	\subcaptionbox{2nd truncation}{
		\includegraphics[height=0.3\linewidth]{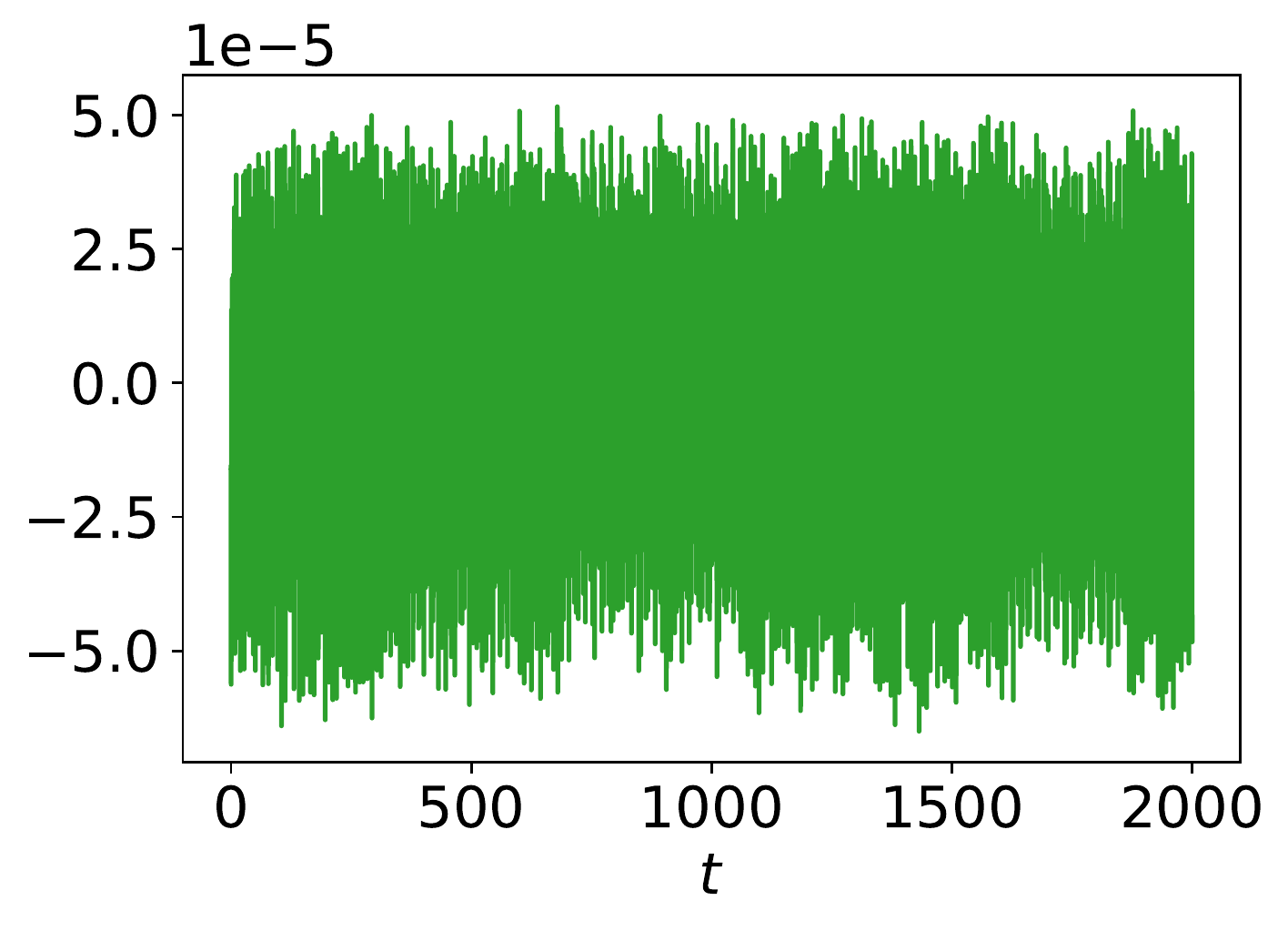}
	}
	
	\caption{H\'enon--Heiles experiment.
		Truncations of the formal power series $\widetilde{\overline{H}}$ are increasingly better conserved along the motion predicted by SSI. This justifies using a system with Hamiltonian $\widetilde{\overline{H}}^{[2]}$ to analyse the dynamical system that describes the numerical motions obtained with SSI.
		(Arithmetic means are subtracted before plotting the energy behaviour.)
	}\label{fig:HenonHeilsHRecovery}
\end{figure}

To analyse the behaviour of numerical motions predicted with SSI and to demonstrate SSI's potential for system identification, we apply the backward error analysis formula \eqref{eq:tildeH} to the learned Hamiltonian $\overline{H}$ to obtain the formal power series $\widetilde{\overline{H}}$.
Truncations of $\widetilde{\overline{H}}$ now describe up to the truncation error the Hamiltonian of a dynamical system which governs the numerical motion. Since the truncation to second order $\widetilde{\overline{H}}^{[2]}$ is sufficiently well preserved for our purposes (\cref{fig:HenonHeilsHRecovery}),
we use $\widetilde{\overline{H}}^{[2]}$ to recover the potential $V$ as $\tilde{\overline{V}}^{[2]}(q) := \widetilde{\overline{H}}^{[2]}(q,(0,0))$. Contour plots of the exact and recovered potential are visually indistinguishable (\cref{fig:HenonHeils}).
Notice that also the separatrix is reproduced correctly, which is related to the favourable energy conservation properties of SSI and explains the excellent long term behaviour for motions initialised close to the separatrix. To quantify the difference between the numerical and exact phase portrait, we consider a uniform distribution on the training domain $M$ and compute the standard deviation $\sigma = \sigma(H-\widetilde{\overline{H}}^{[2]})$ approximated using a uniform 20x20x20x20 grid. We obtain $\sigma < 7 \cdot 10^{-4}$. 
This confirms that the numerical and exact phase portrait are close on $M$, which guarantees qualitatively correct long-term behaviour.
%

\section{Future work}\label{sec:FutureWork}

Symplectic Shadow Integration is compatible with techniques to incorporate symmetries into learning processes.\cite{ridderbusch2021learning} If symmetry groups of the system are known, equivariant kernels adapted to the system's symmetries can be used to model $\overline{H}$.
If the symmetries are symplectic and preserved by the integrator, the numerical flow will share the symmetry group of the exact flow. This is especially valuable when the action of the symmetry group is symplectic. In such cases symmetries relate to conserved quantities by Noether's theorem,\cite{mansfield2010} which are then picked up by the numerical flow. This leads to excellent preservation of topological properties of the phase portrait and can be used, for instance, to conserve complete integrability.\cite{GeomIntegration}

In \cref{sec:learnH} the inverse modified Hamiltonian $\overline{H}$ is obtained as the mean function of a Gaussian process $\hat{\overline{H}}$.
The technique can easily be extended to general kernel based method with continuously differentiable kernels. For the system identification part, higher differentiability may be required.  
When Gaussian Processes are used, information about uncertainty of $\hat{\overline{H}}$ could be used to predict uncertainty of trajectories computed with SSI.
Moreover, SSI can be extended to non-canonical Hamiltonian systems, if a symplectic integration method is known for the system type, which is the case for spin systems,\cite{Robert2014} for instance.
Furthermore, building on variational backward error analysis techniques,\cite{Vermeeren2017} the authors plan to develop a variational version of SSI which learns inverse modified Lagrangian functions that compensate discretisation errors of variational integrators.\\



For a summary of the findings, also see our  poster.\cite{SymplecticShadowIntPoster}

\section*{Data Availability Statement}

The data that support the findings of this study are openly available from the repository Christian-Offen/symplectic-shadow-integration \cite{ShadowIntegratorSoftware} (\url{https://github.com/Christian-Offen/symplectic-shadow-integration},\\
\url{https://dx.doi.org/10.5281/zenodo.5159766}).

\appendix
\section{Relation of inverse modified Hamiltonians to generating functions for the implicit midpoint rule}\label{App:GenFuncMP}

\noindent {\em Lemma}. Consider the application of the implicit midpoint rule to a Hamiltonian system with step size $h$. Denote its inverse modified Hamiltonian by $\overline {H}$. 
Around any point in the phase space, there exists a generating function $S$ of the exact flow map $\phi_h$ such that the scaled inverse modified Hamiltonian $h\overline{H}$ fulfils the same partial differential equations as $S$.

\begin{proof}
	The graph $\Gamma = \{(q,p,\overline{q},\overline{p}) \, | \, (\overline{q},\overline{p})=\phi_h(q,p) \}$ is an embedded Lagrangian submanifold in a symplectic manifold with symplectic structure $\Omega = \sum_{j=1}^{n} (\d \overline{q}^j \wedge \d \overline{p}^j - \d q^j \wedge \d p^j)$. The following 1-form $\alpha$ is a primitive of $\Omega$, i.e.\ $\d \alpha = \Omega$:
	\[
	\alpha = \sum_{j=1}^n \left( (\overline{q}^j-q^j) \d \left( \frac{\overline{p}^j+p^j}{2} \right)  - (\overline{p}^j-p^j) \d \left( \frac{\overline{q}^j+q^j}{2} \right) \right).
	\]
	The (pullback of the) 1-form $\alpha$ is closed on $\Gamma$, since $\Gamma$ is a Lagrangian submanifold. Therefore, primitives $S$ with $\d S = \alpha$ exist locally on $\Gamma$. For small step sizes, $\frac{\overline{q}+q}{2},\frac{\overline{p}+p}{2}$ constitutes a local coordinate system on $\Gamma$. Writing
	\begin{align*}
		\d S &= \sum_{j=1}^n \frac{\p S}{\p \left( \frac{\overline{q}^j+q^j}{2} \right)}\left( \frac{\overline{q}^j+q^j}{2}, \frac{\overline{p}^j+p^j}{2}\right) \d \left( \frac{\overline{q}^j+q^j}{2} \right) \\
		&+\sum_{j=1}^n\frac{\p S}{\p \left( \frac{\overline{p}^j+p^j}{2} \right)}\left( \frac{\overline{q}^j+q^j}{2}, \frac{\overline{p}^j+p^j}{2}\right) \d \left( \frac{\overline{p}^j+p^j}{2} \right)
	\end{align*}
	and comparing the coefficients of the total differentials of $\alpha$ and $\d S$ yield equations, which after substituting $h \overline{H}$ for $S$ read
	\begin{equation}\label{eq:IMHbar}
\begin{split}
	\overline{q} &= q+h {\nabla_p} \overline{H}\left( \frac{\overline{q}^j+q^j}{2}, \frac{\overline{p}^j+p^j}{2}\right)\\
\overline{p} &= p-h {\nabla_q} \overline{H}\left( \frac{\overline{q}^j+q^j}{2}, \frac{\overline{p}^j+p^j}{2}\right).		
\end{split}
	\end{equation}
Here, ${\nabla_q} \overline{H}$ and ${\nabla_p} \overline{H}$ denote derivatives of $\overline{H}$ with respect to the first input argument or second input argument, respectively.
The system \eqref{eq:IMHbar} coincides with the implicit midpoint rule applied to $\overline{f} = J^{-1} \nabla \overline{H}$. 
\end{proof}

\bibliography{bib}

\end{document}

%% file: ORCID_logo.tex
\definecolor{orcidlogocol}{HTML}{A6CE39}
\tikzset{
  orcidlogo/.pic={
    \fill[orcidlogocol] svg{M256,128c0,70.7-57.3,128-128,128C57.3,256,0,198.7,0,128C0,57.3,57.3,0,128,0C198.7,0,256,57.3,256,128z};
    \fill[white] svg{M86.3,186.2H70.9V79.1h15.4v48.4V186.2z}
                 svg{M108.9,79.1h41.6c39.6,0,57,28.3,57,53.6c0,27.5-21.5,53.6-56.8,53.6h-41.8V79.1z M124.3,172.4h24.5c34.9,0,42.9-26.5,42.9-39.7c0-21.5-13.7-39.7-43.7-39.7h-23.7V172.4z}
                 svg{M88.7,56.8c0,5.5-4.5,10.1-10.1,10.1c-5.6,0-10.1-4.6-10.1-10.1c0-5.6,4.5-10.1,10.1-10.1C84.2,46.7,88.7,51.3,88.7,56.8z};
  }
}

\newcommand\orcidicon[1]{\href{https://orcid.org/#1}{\mbox{\scalerel*{
\begin{tikzpicture}[yscale=-1,transform shape]
\pic{orcidlogo};
\end{tikzpicture}
}{|}}}}